\documentclass[12pt]{amsart}

\usepackage{amssymb,amsmath,bbm}
\usepackage[matrix,arrow,graph]{xy}
\numberwithin{equation}{subsection}

\newtheorem{thm}[equation]{Theorem}
\newtheorem{lemma}[equation]{Lemma}
\newtheorem{prop}[equation]{Proposition}
\newtheorem{cor}[equation]{Corollary}
\theoremstyle{definition}

\newtheorem{defn}[equation]{Definition}

\newtheorem{ex}[equation]{Example}

\newtheorem*{rmk}{Remark}

\theoremstyle{remark}

\newcommand{\bbC}{{\mathbb C}}

\newcommand{\bbL}{{\mathbb L}}
\newcommand{\bbR}{{\mathbb R}}

\newcommand{\bbZ}{{\mathbb Z}}

\newcommand{\bbP}{{\mathbb P}}

\newcommand{\Z}{{\mathbb Z}}

\newcommand{\C}{{\mathbb C}}

\newcommand{\R}{{\mathbb R}}

\newcommand{\cF}{{\mathcal F}}

\newcommand{\cO}{{\mathcal O}}
\newcommand{\cP}{{\mathcal P}}
\newcommand{\cL}{{\mathcal L}}
\newcommand{\cM}{{\mathcal M}}
\newcommand{\cN}{{\mathcal N}}

\newcommand{\cA}{{\mathcal A}}
\newcommand{\cB}{{\mathcal B}}
\newcommand{\cI}{{\mathcal I}}
\newcommand{\cC}{{\mathcal C}}
\newcommand{\cE}{{\mathcal E}}

\newcommand{\cT}{{\mathcal T}}
\newcommand{\cK}{{\mathcal K}}
\newcommand{\cH}{{\mathcal H}}

\newcommand{\Sh}{\operatorname{Sh}}
\newcommand{\supp}{\operatorname{Supp}}

\newcommand{\im}{\operatorname{Im}}

\newcommand{\Hom}{\operatorname{Hom}}
\newcommand{\Ext}{\operatorname{Ext}}

\newcommand{\Symm}{\operatorname{Sym}}

\newcommand{\Span}{\operatorname{Span}}

\newcommand{\Sym}{\operatorname{Sym}}

\newcommand{\St}{\operatorname{St}}

\renewcommand{\check}[1]{{#1}^\vee}

\newcommand{\la}{\langle}
\newcommand{\ra}{\rangle}

\newcommand{\wt}{\widetilde}

\newcommand{\ol}{\overline}

\newcommand{\rank}{\mathop{\rm rank}\nolimits}
\renewcommand{\ker}{\mathop{\rm ker}\nolimits}
\newcommand{\coker}{\mathop{\rm coker}\nolimits}
\newcommand{\id}{{\rm id}}

\renewcommand{\hom}{\mathop{\rm hom}\nolimits}
\newcommand{\End}{\mathop{\rm End}\nolimits}

\newcommand{\PP}{{\mathbb P}}

\newcommand{\cal}{\mathcal}

\newcommand{\cZ}{{\cal Z}}

\newcommand{\bdy}{{\partial}}
\newcommand{\sig}{{\sigma}}
\newcommand{\Sig}{{\Sigma}}
\newcommand{\udot}{{\scriptscriptstyle \bullet}}

\newcommand{\Gr}{\mathop{\rm Gr}\nolimits}
\newcommand{\ext}{\mathop{\rm ext}\nolimits}
\newcommand{\co}{\text{\rm co-}}
\renewcommand{\mod}{\text{\rm -mod}}
\newcommand{\Mod}{\text{\rm -Mod}}
\newcommand{\mof}{\text{\rm -mod}_{\text{\it f}}}
\newcommand{\mcf}{\text{\rm -mod}_{\text{\it cf}}}

\newcommand{\Mof}{\text{\rm -Mod}_{\text{\it f}}}
\newcommand{\Mcf}{\text{\rm -Mod}_{\text{\it cf}}}

\newcommand{\Lie}{\mathop{\rm Lie}}
\newcommand{\Pure}{\mathop{\rm Pure}}

\newcommand{\bC}{{\mathbf C}}
\newcommand{\bD}{{\mathbf D}}

\newcommand{\bM}{{\mathbf M}}

\newcommand{\bR}{{\mathbf R}}
\newcommand{\bS}{{\mathbf S}}

\newcommand{\tmop}[1]{\operatorname{#1}}

\newcommand{\eend}{\tmop{end}}

\newcommand{\m}{\mathfrak m}
\newcommand{\For}{\mathop{\mathit For}}
\newcommand{\lb}{\{}
\newcommand{\rb}{\}}

\begin{document}
 
\title[Equivariant-constructible Koszul duality]
{Equivariant-constructible Koszul duality for dual toric varieties}
\date{\today}
\author{Tom Braden}
\address{Dept.\ of Mathematics and Statistics\\
         University of Massachusetts, Amherst}
\email{braden@math.umass.edu}
\author{Valery A.~Lunts}
\address{Department of Mathematics, Indiana University,
Bloomington, IN 47405, USA}
\email{vlunts@indiana.edu}

\thanks{The first named author was partially supported by NSF grant
 DMS-0201823.\\
The second named author was partially supported by NSA grant MDA904-01-1-0020
and CRDF grant RM1-2405-MO-02}

\begin{abstract} For affine toric varieties $X$ and $\check X$ 
defined by dual cones, we define an equivalence of categories between
mixed versions of the equivariant derived category $D^b_T(X)$ and 
the derived category of sheaves on $\check X$ which are locally constant
with unipotent monodromy on each orbit.  This equivalence satisfies the
Koszul duality formalism of Beilinson, Ginzburg, and Soergel.
\end{abstract}

\subjclass{14M25; 16S37, 55N33, 18F20}
% Toric Varieties, Koszul algebras, IC, Sheaves
\maketitle
\bibliographystyle{amsabbrv}
 
\setcounter{tocdepth}{1}
\tableofcontents

\section{Introduction}

\subsection{}
Let $T$ and $\check T$ be dual complex tori. The ring $\cO (T)$ of
regular functions on $T$ is canonically isomorphic to the group ring
$\bbC [\pi _1(\check T)]$ of the fundamental group $\pi _1(\check T)$.
Thus the category of quasi-coherent sheaves on $T$ is identified with
the category of local systems on $\check T$. Let $\mathfrak{t}$ be the
Lie algebra of $T$. The (evenly graded) algebra $A$ of polynomial
functions on $\mathfrak{t}$ is canonically isomorphic to the
equivariant cohomology ring $H^*_T(pt)$.
The exponential map
$$exp: \mathfrak{t}\to T$$ identifies 
$A$-modules supported at the origin with quasi-coherent sheaves on 
$T$ supported at the identity and hence with unipotent local systems 
on $\check T$.

Notice that under this correspondence the logarithm of
the monodromy operator becomes the multiplication by the first Chern
class of a line bundle. This phenomenon also occurs in 
mirror symmetry, where the logarithm of the monodromy of 
the Gauss-Manin local system around certain loops in the moduli space
of complex structures of a Calabi-Yau manifold is identified with
multiplication by a Chern class on the cohomology of  
the mirror Calabi-Yau.  This is worked out in more generality for
the case of complete intersections in toric varieties in
\cite{Ho}.   
%So it is proper to consider our eventual correspondence
%between equivariant sheaves on a toric variety $X$ and (unipotent)
%locally constant sheaves on the dual toric variety $\check X$ as
%having two components: mirror symmetry and Koszul duality....

\subsection{} \label{intro to main results}
In this paper we extend the above idea to affine toric varieties.  
Let $X$ be an affine $T$-toric variety whose fan
consists of a cone $\sig$ and all its faces, 
and let $\check X$ be the toric variety whose fan 
consists of the dual cone $\check \sig$ and all its faces.
The $T$-orbits of $X$ are indexed by the faces $\tau$ of the
cone $\sig$.  The $T$-orbit $O_\tau \subset X$ is 
identified with a quotient $T/T_\tau$ by a subtorus $T_\tau$.
The dual variety has a corresponding orbit $O_{\tau^\bot}$ which is
isomorphic as a $\check T$-space to $\check T/T_\tau^\bot$,
where $T_\tau^\bot$ is the  ``perpendicular'' subtorus
whose Lie algebra is the annihilator of 
$\mathfrak{t}_\tau = \Lie T_\tau$.

The ring of regular functions on $T_\tau$ is canonically
isomorphic to the group ring of $\pi_1(\check T/T_\tau^\bot)$.
Let $A_\tau$ be the algebra of polynomial functions on 
$\mathfrak{t}_\tau$; it is canonically  isomorphic to
the equivariant cohomology $H^*_T(T/T_\tau)$.
In the same way as before, the exponential map identifies 
$A_\tau$-modules supported at the origin with unipotent 
local systems on $O_{\tau^\bot}$. 

We use this idea (together with a more combinatorial
duality, see \S\ref{combinatorial Koszul}) to relate 
$T$-equivariant sheaves on $X$ to complexes of sheaves 
with unipotent monodromy on $\check X$.  We define an 
equivalence of triangulated categories $K$, which fits 
into the following diagram of categories and functors. 
\begin{equation} \label{main diagram}
\xymatrix{
D^b(\cA\mof) \ar[r]^K_\sim\ar[d]^{F_T} & D^b(LC_\cF(X^\vee))\ar[d]^{F_{cf}} \\
 D^b_T(X) & D^b(LC_{cf}(X^\vee))
}\end{equation}
The categories on the bottom are topological categories of sheaves
on the dual toric varieties $X$ and $X^\vee$.  The categories 
above them are ``mixed'' versions of these categories, where 
objects have been given an (extra) grading, 
and the vertical functors forget the grading. 
On the left hand side we have $T$-equivariant sheaves on $X$,
while on the right we have complexes of orbit-constructible sheaves
on $\check X$ with unipotent monodromy. 

Let us describe these categories in more detail.
$D^b_T(X)$ is the (bounded, constructible) 
$T$-equivariant derived category of sheaves on
$X$ defined in \cite{BL}.  By results of \cite{L} it is equivalent
to a full subcategory of 
the category of differential graded modules (DG-modules) over
a sheaf $\cA = \cA_{[\sig]}$ of rings on the finite poset
$[\sig]$ of faces of $\sig$.  Sections of this sheaf on a 
face $\tau\prec \sig$ are complex valued polynomial functions
on $\tau$, graded so linear functions have degree $2$. 

Our mixed version of this category is $D^b(\cA\mof)$, the
derived category of finitely generated graded $\cA$-modules; it 
has two gradings -- the module grading, and the grading in the complex.
The forgetful functor $F_T$ combines these two gradings into 
the single grading on DG-modules.  There is a ``twist'' 
automorphism of $D^b(\cA\mof)$, denoted $\la 1\ra$,
which shifts both the complex and module gradings so that
$F_T\la 1\ra = F_T$.

The right side of \eqref{main diagram} involves sheaves on 
$\check X$. $LC_{cf}(\check X)$
denotes the category of sheaves of $\C$-vector spaces 
on $\check X$ which are
\begin{enumerate}
\item[(1)] ``locally constant'' -- the restriction to each 
$T$-orbit $O \subset \check X$ is a local system, 
\item[(2)] ``unipotent'' -- for any orbit $O$ and any 
$\gamma\in \pi_1(O)$, the action of $\gamma-1$ on
the stalk at a point $p \in O$ is locally nilpotent.
\item[(3)] ``cofinite'' -- the $\pi_1(O)$ invariants
of the stalk at $p$ is finite dimensional, for all 
orbits $O$ and $p\in O$.
\end{enumerate} 
$D^b(LC_{cf}(\check X))$ is then the derived category of
this abelian category.  Note that although 
objects of $D^b(LC_{cf}(\check X))$ are locally constant on orbits,
they are not constructible in the usual sense, since the stalks
need not be finite-dimensional.  As we will see, though, conditions 
(1)--(3) imply that these objects are
still well-behaved.  In particular, the full subcategory of 
objects all of whose stalk cohomology groups are 
finite-dimensional is equivalent to a full subcategory of 
the usual constructible derived category. 
See Proposition \ref{constructible = finite length}.
 
To define the mixed version $LC_\cF(\check X)$ of 
$LC_{cf}(X^\vee)$, we use a self-map $\cF\colon
\check X \to \check X$ which is a lift of the Frobenius
map to characteristic zero defined for toric varieties.
An object in $LC_\cF(\check X)$ is an object $S \in LC_{cf}(\check X)$
together with an isomorphism $\theta \colon \cF ^{-1}F\to F$
whose eigenvalues on the stalks at points of $(\check X)^\cF$ 
are powers of $2^{1/2}$.  We again have a ``shift of grading''
functor $\la 1\ra$, which multiplies $\theta$ by $2^{1/2}$.
It is clear that $F_{cf}\la 1 \ra = F_{cf}$. 

\subsection{$t$-structures} 
All four categories in \eqref{main diagram} 
come equipped with natural perverse
$t$-structures, whose abelian cores we denote by $P_T(X)$, $P(\cA)$,
$P_{cf}(\check X)$, and $P_\cF(\check X)$.  These $t$-structures are
particularly nice, in that each triangulated category is equivalent to
the bounded derived category of its core: $D^b_T(X) \cong
D^b(P_T(X))$, etc.
%The $t$-structure on $D^b_T(X)$ was defined in \cite{BL}, and
%the $T$-structures on $D^b(LC_\cF(X^\vee))$ and $D^b(LC_{cf}(X^\vee))$
%are defined in the usual way using degree restrictions on 
%stalks and costalks.  The $t$-structure on $D^b(\cA\mof)$ can 
The forgetful functors $F_T$, $F_{cf}$ are $t$-exact, so 
they restrict to functors $P(\cA)\to P_T(X)$ and 
$P_\cF(\check X) \to P_{cf}(\check X)$.

%In fact $F_T$ and $F_{cf}$ are
%``strongly $t$-exact'' -- meaning that $\cM^\udot \in D^b(\cA\mof)$ is
%perverse if and only if $F_T(\cM^\udot)$ is perverse in $D^b_T(X)$, 
%and similarly for $F_{cf}$.

The twist functors $\la 1\ra$ are also $t$-exact, and so they give
automorphisms of $P(\cA)$ and $P_\cF(X^\vee)$.  This induces
bijections between isomorphism classes of ``ungraded'' simple objects
and ``graded'' simples up to twists:
\begin{eqnarray*} 
Irr(P_T(X)) & \leftrightarrow & Irr(P(\cA))/\Z,\\
Irr(P_{cf}(\check X)) & \leftrightarrow & Irr(P_\cF(\check X))/\Z.
\end{eqnarray*}
The set $Irr(P_T(X))$ consists of all equivariant IC-sheaves supported
on the closures of $T$-orbits. For each IC-sheaf $IC_T(\ol{O_\tau})$
we will fix a certain lift to $P(\cA)$, which we denote $\cL^\tau$;
it is a complex of $\cA$-modules with nonzero cohomology in a single degree.
Up to a grading shift, it is the combinatorial equivariant
intersection cohomology sheaf studied in \cite{BBFK, BrLu, Ka}.

%; we call these simples ``pure of
%weight $0$''. 

\subsection{} \label{injectives and projectives}
The abelian category $P(\cA)$ has enough projective objects.  Let $\cP\in
P(\cA)$ be a projective cover of $\oplus_{\tau\in [\sig]}\cL^\tau$,
and let $R$ be the opposite ring to the graded ring
\[\eend(\cP) := \oplus_{i\ge 0} \Hom_{P(\cA)}(\cP,\cP\la i\ra).\]  
In the usual way we see that $P(\cA)$ is equivalent to $R\mof$, the
category of finitely generated $R$-modules.  Furthermore, $P_T(X)$ is
equivalent to $R\Mof$, the category of finitely generated ungraded
$R$-modules.  It follows that we have equivalences $D^b(\cA\mof)\cong
D^b(R\mof)$, $D^b_T(X)\cong D^b(R\Mof)$.  With respect to these
equivalences, $F_T$ is the functor of forgetting the grading.

A similar story holds on the right-hand side of (\ref{main diagram}).
The simple objects in $P_{cf}(\check X)$ are the intersection
cohomology complexes $IC^\udot(\ol{O_\alpha})$ for $O_\alpha \subset
\check X$ a $T$-orbit.  We will single out a distinguished 
lift $L^\udot_\alpha\in P_\cF(\check X)$ of $IC^\udot(\ol{O_\alpha})$.
%, which we call ``pure of weight 0''.

There are enough injectives in $P_\cF(\check X)$; let $I$ denote the
injective hull of $\oplus L^\udot_\alpha$, and put $R^\vee =
\eend(I)^{opp}$.  The functor $\oplus_{i\ge 0} \Hom_{P_\cF(\check
  X)}(-,I\la -i\ra)^*$ gives an equivalence between $P_\cF(\check X)$
and $\check R\mcf$, the category of ``co-finite'' graded $\check
R$-modules (see \S\ref{module conventions}).  
%The ungraded version 
%$F_{cf}(I)$ of $I$ is again injective, and so 
Similarly $P_{cf}(\check X)$ is equivalent to the
category $\check R\Mcf$ of ungraded co-finite modules.  With respect
to these equivalences, $F_{cf}$ is the functor of forgetting the
grading.

The full subcategory $P_{u,fl}(\check X)$ of $P_{cf}(\check X)$
consisting of objects of finite length is the full subcategory of the
usual topological category of orbit-constructible perverse sheaves
consisting of objects with unipotent monodromy.  Further, the category
$P_{\cF,fl}(\check X)$ of finite length objects in $P_\cF(\check X)$
is a mixed version of $P_{u,fl}(\check X)$.  These categories are
equivalent to finite-dimensional ungraded and graded $\check
R$-modules, respectively.

\subsection{} 
The functor $K$ which relates the two sides of 
(\ref{main diagram}) is {\em not\/} $t$-exact, but it 
does have an interesting relationship with the 
$t$-structures: it is a Koszul equivalence. 
Roughly this means that $K$ takes simples of weight $0$ in 
$P(\cA)$ to indecomposable injectives in $P_\cF(\check X)$ 
and indecomposable projectives in $P(\cA)$ to simples
in $P_\cF(\check X)$.  It also implies that  
$R$ and $\check R$ are Koszul graded rings, and are naturally 
Koszul dual to each other, in the sense of \cite{BGS}.
We explain this in more detail in \S\ref{def of Koszul functor} below.

\subsection{}  A similar Koszul functor was constructed 
for the varieties $X$ and $\check X$ in \cite{Br}.  In that 
construction the source and target categories were 
combinatorially defined triangulated categories $\bD_\Phi(X)$
and $\bD_{\check \Phi}(\check X)$, which model mixed sheaves
on $X$ and $X^\vee$ with ``conditions at infinity'' described
by auxiliary data $\Phi$, $\Phi^\vee$.  This auxiliary choice
(essentially the choice of a toric normal slice to each stratum)
is somewhat artificial, and as a result it is not clear how
to relate these categories directly to a topological category,
although the corresponding abelian category of perverse 
objects in $\bD_\Phi(X)$ is a mixed version of a category 
$\cP_\Phi(X)$ of perverse sheaves on $X$.  

The construction in this paper
removes these defects, and it is our hope that this more
canonical approach will in turn inspire a more intrinsic point
of view on this phenomenon, in which the duality functor is
defined directly, perhaps in terms of filtered $D$-modules.

\subsection{Ideas for future work}
We expect that our work can be extended in several directions.

a) We hope that a particular instance of our constructible-equivariant
correspondence can be considered as a ``limit case'' of the mirror
symmetry between dual families of Calabi-Yau hypersurfaces in dual toric 
varieties.

b) We believe that the same constructible-equivariant component should be
present in the Koszul duality on flag manifolds constructed in
[BGS]. This should be in agreement with Soergel's conjectures [S].

c) One should be able to extend to toric varieties the full
correspondence (local systems on
$T$) $\leftrightarrow$ (quasi-coherent sheaves on $\check T$). In
our work we restricted ourselves to unipotent local systems and
quasi-coherent sheaves supported at the identity. The language of
configuration schemes [L2] may be appropriate in this problem.
%...

\subsection{} We briefly describe the structure of the paper.
Section 2 contains some basic background from homological
algebra, including definitions on graded modules, 
discussion of Koszul equivalences, and mixed categories and 
gradings.  Section 3 introduces the formalism of sheaves on
fans considered as finite partially ordered sets, and defines
three sheaves of rings on fans which are important later. 

In Section 4 we consider sheaves on toric varieties which are
locally constant on orbits.  We prove that the derived category
of these sheaves is the same as the category of complexes of
sheaves with locally constant cohomology; this means that locally 
constant sheaves have enough flexibility for our homological 
calculations.  We also show that the category $LC(X)$ of 
locally constant sheaves is equivalent to comodules over a sheaf 
of rings $\cB$ defined in section 3.

In section 5 we define our ``mixed'' version 
$LC_\cF(X)$ of locally constant sheaves, and show that they are equivalent 
to graded comodules over $\cB$.  We define a perverse $t$-structure 
on $D^b(LC_\cF(X))$ and prove some basic properties of perverse 
objects, including the local purity of simple objects.  This
purity allows us to define a mixed structure on the category
of perverse objects, which is the first step to proving that
they are equivalent to modules over some graded algebra.

Section 6 turns to the equivariant side of our picture.  We first
describe a topological realization functor from
complexes of $\cA$-modules on a fan $\Sig$ to equivariant complexes on 
the corresponding toric variety $X$; this was originally defined in \cite{L}.  
We next study the homological algebra of $\cA$-modules; the main result is
that $D^b(\cA\mof)$ is equivalent to the 
homotopy category of complexes of ``pure'' $\cA$-modules, which 
are direct sums of shifts of the combinatorial
intersection cohomology sheaves $\cL^\tau$ studied in \cite{BBFK,BrLu}. 

In Section 7 we finally define our toric Koszul functor $K$
and prove that it has the asserted properties.  In fact 
the existence of enough projectives in $P(\cA)$ and 
enough injectives in $P(LC_\cF(\check X))$ is deduced
by applying $K$ and $K^{-1}$ to the appropriate simple
perverse objects, and all the assertions of 
\S\ref{injectives and projectives} are proved here.

We banish a few technical proofs to Section 8.

\subsection{Acknowledgments}
The first author would like to thank David Cox for suggesting
that the results of \cite{Br} should be reformulated equivariantly.
The second author would like to thank V.\ Golyshev for stimulating
discussions.

\section{Ideas from homological algebra}
\subsection{Conventions on graded rings and modules} \label{module conventions}

Fix a field $k$, and let $R = \oplus_{n\ge 0} R_n$ be a positively
graded $k$-algebra whose zeroth graded piece $R_0$ is isomorphic to
$k^{\oplus s}$ for some $s \ge 1$.  Suppose that all graded pieces
$R_n$ are finite-dimensional (this holds if $R$ is either left or
right Noetherian, for instance).

Let $R\mod$, $R\Mod$ denote the abelian categories of graded (resp.\ 
ungraded) left $R$-modules.  Let $R\mof$, $R\Mof$ be their respective
full subcategories of finitely generated modules; they are abelian
subcategories if and only if $R$ is left Noetherian.

The shift of grading functors $\la j\ra, j\in \Z$ act on $R\mod$ by
$(M\la j\ra)_n = M_{n+j}$ (note that this is the opposite convention
to \cite{BGS}).  They preserve the subcategory $R\mof$.

Given $M$ in $R\mod$, define the ``graded dual'' $M^*$ of $M$ by
$(M^*)_n = \Hom_k(M_{-n}, k)$.  Then $M\mapsto M^*$ is a functor
$R\mod \to (R^{opp})\mod^{opp}$.  Put $R^\circledast := (R^{opp})^*$;
it is an injective object of $R\mod$.

We will also need ``cofinite'' graded $R$-modules, the dual notion to
finitely generated modules.  If $n \ge 0$, put $R_{> n} = \oplus_{k>
  n}R_k$.
\begin{prop} Let $M \in R\mod$ or $R\Mod$.  The following are equivalent:
\begin{enumerate}
\item $\dim_k \{m\in M \mid R_{>0}\cdot m = 0\}<\infty$, and every
  $m\in M$ is annihilated by some $R_{> n}$.
\item $M^*$ is contained in a finite direct sum of shifted copies of
  $R^{\circledast}$.
\end{enumerate}
\end{prop}
Such modules are called $R_{>0}$-cofinite, or simply ``cofinite''.
Let $R\mcf$, $R\Mcf$ denote the category of cofinite graded and ungraded 
$R$-modules, respectively.  
The graded dual gives an equivalence of categories $R\mof
\to (R^{opp})\mcf^{opp}$.  It follows that $R\mof$ is an abelian
subcategory of $R\mod$ if and only if $R$ is right Noetherian.

\subsection{Koszul functors and Koszul duality} \label{def of Koszul functor}
We present the ideas of Koszul duality on derived categories at a
level of generality appropriate for our purposes.  For a more general
discussion, see \cite{BGS}.

Fix a field $k$. Let $R$ and $R^\vee$ be algebras of the type
considered in the previous section.

\begin{defn} \label{Koszul functor definition}
  A covariant functor
\[K\colon  D^b(R\mof) \to D^b(R^\vee\mcf)\]
is a {\em Koszul equivalence} if the following are satisfied:
\begin{enumerate}
\item $K$ is a triangulated equivalence of categories.  In particular
  $K(M[1]) = (KM)[1]$ for all $M\in D^b(R\mof)$.
\item For all $M \in D^b(R\mof)$, we have
\[K(M\la 1\ra) = (KM)\la -1 \ra[1].\]
\item $KR_0 \cong (R^\vee)^\circledast$.
\item $KR \cong (R^\vee_0)^\circledast$.
\end{enumerate}
\end{defn}

Since $R_0$, $R^\vee_0$ are semisimple, $R$ is projective and \ 
$(R^\vee)^\circledast$ is injective, the conditions (3) and (4) can be
replaced by the following.
\begin{enumerate}
\item[($3'$)] $K$ sends simple objects of grading degree $0$ to
  injective hulls of simples of degree $0$.
\item[($4'$)] $K$ sends projective covers of simples of grading degree
  $0$ to simples of degree $0$.
\end{enumerate}
Here we use the standard embeddings of $R\mof$, $\check R\mcf$ into
their derived categories as complexes with cohomology only in degree
$0$.

\begin{thm} \label{Koszul functors and rings}
If a Koszul functor $K$ exists as in Definition \ref{Koszul functor definition},
  then
\begin{enumerate}
\item[(a)] $R$, $R^\vee$ are Koszul graded rings, i.e.\ 
\[\Ext^i_{R}(R_0,R_0\la j\ra) = 0\;\text{for}\; i\ne -j,\]
and similarly for $R^\vee$.
\item[(b)] $R^\vee$ is the Koszul dual ring to $R$, i.e.\ we have an
  isomorphism of rings
\[\label{Koszul dual ring}
R^\vee \cong \oplus_{i\ge 0} \Ext^i_R(R_0,R_0\la -i\ra).\]
\end{enumerate}
\end{thm}

The proof is immediate from the definition \ref{Koszul functor
  definition}.

The conclusion (a) implies in particular that $R$ is a quadratic algebra, with
generators in degree $1$ and relations in degree $2$.  (b)
implies that $R$ and $R^\vee$ are quadratic dual rings: $R_0 =
R^\vee_0$ canonically, $R_1$ and $R_1^\vee$ are dual $R_0$-modules,
and the relations for $R$ and $R^\vee$ are orthogonal.  See \cite{BGS}
for more precise statements and a proof.

\begin{rmk}
  The original example of a Koszul equivalence was defined by
  Bernstein, Gelfand, and Gelfand \cite{BGG} for $R$ a polynomial ring
  and $R^\vee$ the dual exterior algebra.  This example is related to
  a duality between equivariant and ordinary cohomology, and underlies
  the ``local'', one-orbit case of our toric Koszul duality.
  
  In \cite{BGS} it is shown that under mild finiteness conditions
  (e.g.\ if $\dim_k R^\vee <\infty$, so $R^\vee\mcf = R^\vee\mof$),
  then a Koszul dual pair of rings $(R, R^\vee)$ gives rise to a
  Koszul equivalence $D^b(R\mof)\to D^b(R^\vee\mof)$.  We are taking
  the opposite point of view and considering the functor $K$ as 
  the primary object.
\end{rmk}

\subsection{Mixed categories} \label{mixed categories}
We will need the notion of a ``mixed'' abelian category.  This
generalizes the category of finitely generated graded modules over a
finitely generated positively graded ring.  We mostly follow
\cite{BGS}, but we do not wish to assume our algebras are
finite-dimensional, so our abelian categories are not assumed to be
Artinian.

Fix a field $k$.  Consider triples $(\bM, W_\udot, \la 1\ra)$, where
\begin{itemize}
\item $\bM$ is an abelian $k$-category. % with finite dimensional $\Hom$s,
\item $\la 1\ra$ is an automorphism of $\bM$, and
\item For each $M \in \bM$, $\{W_jM\}_{j\in \Z}$ is a functorial
  increasing filtration of $M$.
\end{itemize}
We call such a triple a {\em mixed} category if the following are
satisfied:
\begin{enumerate}
\item The filtration $W$ is strictly compatible with morphisms, so
  $\Gr^W_j = W_j/W_{j-1}$ is an exact functor.
\item For any $M \in \bM$, we have $W_j(M\la 1\ra) = W_{j-1}(M)\la
  1\ra$.
\item If $\Gr^W_jM = 0$ for $j\ne w$ (we call such an object {\em pure
    of weight $w$}), then $M$ is a finite direct sum of simple
  objects.
\item There are only finitely many isomorphism classes of
simples of weight $0$. %(this assumption is only included for 
%simplicity).
\end{enumerate}
Define automorphisms $\la n\ra$, $n\in \Z$ of $\bM$ by taking powers:
$\la n\ra = \la 1\ra^n$.

We say an object $M$ of $\bM$ has weights $\le j$ (resp.\ has weights
$\ge j$) if $W_jM = M$, (resp.\ $W_{j-1}M = 0$).  If both hold, we say
$M$ is pure of weight $j$; such an object is semisimple of finite
length by (3).

%As another consequence of the
%above axioms, we see that if $X \in \bM$ has weights
%$\le k$ and $Y \in \bM$ has weights $\ge {k+1}$,
%then $\Hom_{\bM}(X,Y)$, $\Hom_{\bM}(Y,X)$, and
%$\Ext^1_{\bM}(X,Y\la 1\ra)$ all vanish. 

Given a mixed category $(\bM,W_\udot,\la 1\ra)$, and objects $X,Y \in
\bM$, define the graded hom and graded ext by
\begin{align*}
  \hom(X,Y)_n & = \Hom_\bM(X,Y\la n\ra) \\
  \ext^i(X,Y)_n & = \Ext^i_\bM(X,Y\la n\ra).
\end{align*}
The graded vector space $\eend(X):= \hom(X,X)$ naturally has the 
structure of a graded ring.

Let $L \in \bM$ be the direct sum of one object from each
isomorphism class of weight $0$ simples.  We call a projective object
$M\in \bM$ a {\em mixed projective generator (resp.\ a mixed injective
  generator)} if
\begin{enumerate}
\item $M$ is projective (resp.\ injective),
\item $M/W_{-1}M \cong L$ (resp.\ $W_0M \cong L$), and
\item for any $X \in \bM$ there exist $r_k\in \Z$ and a surjection
  $\oplus_{k=1}^n M\la r_k\ra \to X$ (resp.\ an injection $X \to
  \oplus_{k=1}^n M\la r_k\ra$).
\end{enumerate}  If $P$ is a mixed projective generator, 
then it is is a projective cover of $L$, and the graded endomorphism
ring $\eend(P)$ is positively graded.  If in addition the endomorphisms
of simple objects in $\bM$ are reduced to scalars, then 
$\eend(P)_0$ is isomorphic to $k^r$, where $r$ is the number of
isomorphism classes of weight $0$ simples in $\bM$.  Similar
statements hold for mixed injective generators.

The main examples of mixed categories are categories of graded modules
over graded rings.
%We will call a graded ring $R = \oplus_{j\ge 0} R_j$
%a {\em finitely presented} $k$-algebra if 
%(1) $R_0$ is a finite dimensional semisimple $k$-algebra,
%and (2) $R$ is finitely generated over
%$R_0$.  
Given a positively graded ring $R$ with $R_0$ semisimple, we have a
mixed category $(\bM, W_\udot, \la 1\ra)$ where $\bM$ consists of
graded $R$-modules $M$ with $\dim_k M_j < \infty$ for all $j$, 
$\la 1\ra$ is the degree shift as defined previously, and 
$W_jM = \oplus_{i\ge -j} M_i$.  If $R$
is left (resp.\ right) Noetherian, then this restricts to a mixed
structure on $R\mof$ (resp.\ $R\mcf$).

We want sufficient conditions for a mixed category to be of the form
$R\mof$ or $R\mcf$.
\begin{prop} \label{mixed categories are modules}
  Let $(\bM,W_\udot,\la 1\ra)$ be a mixed category.
\begin{itemize}
\item[(a)] If $\bM$ has a mixed projective generator $P$ and
  $\eend(P)$ is Noetherian, then $\hom(P,-)$ defines an equivalence of
  categories $\bM\to R\mof$, where $R = \eend(P)^{opp}$.
\item[(b)] If $\bM$ has a mixed injective generator $I$ and $\eend(I)$
  is Noetherian, then $\hom(-,I)^*$ defines an equivalence of
  categories $\bM \to R\mcf$, where $R = \eend(I)^{opp}$.
\end{itemize}
In either case the mixed structure on $\bM$ agrees with the one on
graded modules.
\end{prop}

\subsubsection{Gradings on abelian categories} \label{grading section}
Let $\bC$ be an abelian category.  By a {\em pre-grading} on $\bC$ we
mean a collection $(\bM,W_\udot,\la 1\ra, v, \epsilon)$, where
$(\bM,W_\udot,\la 1\ra)$ is a mixed category, $v\colon \bM \to \bC$ 
is an exact functor,
and $\epsilon$ is a natural isomorphism $v \to v\circ \la 1\ra$,
satisfying: (1) $v$ sends simples to simples and (2)
for any $X,Y$ in $\bM$, the map \[\hom_\bM(X,Y) \to
\Hom_\bC(vX,vY)\] induced by $v, \epsilon$ is bijective.

\begin{prop} \label{ungraded modules}
  Let $(\bM,W_\udot,\la 1\ra, v, \epsilon)$ be a pre-grading on $\bC$.
\begin{enumerate}
\item[(a$'$)] If part (a) of Proposition \ref{mixed categories are
    modules} holds, and in addition $vP$ is a projective generator of
  $\bC$, then $\bC$ is equivalent to $R^{opp}\Mof$.
\item[(b$'$)] If part (b) of Proposition \ref{mixed categories are
    modules} holds, and in addition $vI$ is an injective generator of
  $\bC$, then $\bC$ is equivalent to $R^{opp}\Mcf$.
\end{enumerate}
in either case $v$ is the functor of forgetting the grading.
\end{prop}
In either situation (a$'$) or (b$'$) it follows that for any $X,Y$ in
$\bM$ and $i\ge 0$, the induced map
\[\ext^i_\bM(X,Y) \to \Ext^i_\bC(X,Y)\]
is bijective.  Thus our pre-grading is what in \cite{BGS} was termed a
``grading'' on $\bC$.

\subsection{Triangulated gradings} \label{triangulated gradings}
Let $D$ be a triangulated category.  A triangulated grading
on $D$ is defined to be a tuple $(D_m, \la 1\ra, v,\epsilon)$, where
$D_m$ is a triangulated category, $\la 1\ra$ is a triangulated automorphism
of $D_m$, $v\colon D_m\to D$ is a triangulated functor, and
$\epsilon\colon v \to v \circ \la 1\ra$ is a natural isomorphism,
subject to the condition that the induced map
\[\hom_{D_m}(X,Y) \to \Hom_D(vX,vY)\]
is an isomorphism for any $X,Y\in D_m$, where
$\hom_{D_m}$ is defined as in the previous section.

If we have a grading on an abelian category as in the previous 
section, we get a triangulated grading by letting
$D_m = D^b(\bM)$, $D = D^b(\bC)$.  One can also go in
the other direction, starting from a triangulated grading
as defined above, and endowing $D_m$ and $D$
with $t$-structures for which $\la 1\ra$ and $v$ are $t$-exact.  
Letting $\bM$ and $\bC$ be the abelian cores of $D_m$ and $D$, 
respectively, we get functors $\la 1 \ra\colon \bM\to \bM$ and 
$v\colon \bM\to \bC$ as above.  Endowing $\bM$ with a suitable 
mixed structure, we get a pre-grading. If Proposition
\ref{ungraded modules} applies, it is a grading.  In 
\S\ref{main proofs}) we use this approach to 
prove the fact stated in the introduction
that the functors $F_T$ and $F_{cf}$ are gradings on the
appropriate perverse abelian categories.

\section{Sheaves of rings associated to toric varieties}
%Notational changes: 
%  opposite poset is \Phi^o or ^\circ.
%  

\subsection{Ringed quivers}
Let $\Gamma$ be a finite partially ordered set which we consider as a
category: for any $\alpha ,\beta \in \Gamma$ the set of morphisms
$\Hom(\beta ,\alpha)$ contains a single element if $\beta \geq \alpha$
and is empty otherwise.  A (covariant) functor from $\Gamma $ to the
category of rings is called a sheaf of rings on $\Gamma$. Let $\cA=\cA
_{\Gamma}$ be such a sheaf of rings, i.e. $\cA$ is a collection of
rings $\{\cA _{\alpha}\}_{\alpha \in \Gamma}$ with ring homomorphisms
$\phi _{\beta \alpha}\colon \cA _{\beta}\to \cA _{\alpha}$, if $\beta
\geq \alpha$, satisfying $\phi _{\gamma \beta}\phi _{\beta
  \alpha}=\phi _ {\gamma \alpha}$ if $\gamma \geq \beta \geq \alpha$.
We call the pair $(\Gamma ,\cA _{\Gamma})$ a {\it ringed quiver}.

Assume that for every $\alpha \in \Gamma$ there is given a $\cA
_{\alpha}$-module $\cM_{\alpha}$ with a morphism $\psi _{\beta
  \alpha}\colon \cM_{\beta}\to \cM_{\alpha}$ of $\cA _{\beta}$ modules
(for $\beta \geq \alpha$), such that $\psi _{\gamma \beta}\psi _{\beta
  \alpha}=\psi _{\gamma \alpha}$ if $\gamma \geq \beta \geq \alpha$.
This data will be called an $\cA $-module.  If each $\cM_\alpha$ is
finitely generated over $\cA_\alpha$, we call the resulting
$\cA$-module locally finitely generated.  
$\cA$-modules (resp.\ locally finitely generated $\cA$-modules) form
an abelian category which we denote $\cA \Mod$ (resp. $\cA\Mod_f$).
%Let
%$D(\cA \Mod)$ (resp. $D^b(\cA \Mod)$) denote the derived (resp. bounded derived)
%category of complexes in $\cA \Mod}$.

Notice that $\Gamma$ can be viewed as a topological space where $\beta
$ is in the closure of $\alpha $ iff $\beta \geq \alpha$. So the
subsets $[\beta]= \{ \alpha \mid \alpha \leq \beta\}$ are the
irreducible open subsets in $\Gamma$. Then $\cA$ induces a sheaf of
rings on this topological space, so that the corresponding category of
sheaves of modules is equivalent to $\cA \Mod$.

We call $(\Gamma ,\cA )$ a graded ringed quiver if rings $\cA
_{\alpha}$ are graded and $\phi _{\beta \alpha}$ are morphisms of
graded rings.  In this case let $\cA \mod$ (resp. $\cA\mof$) denote
the abelian category of graded $\cA$-modules (resp. locally finitely
generated graded $\cA$-modules) with morphisms of degree zero.

\begin{rmk} The category $\cA\mod$ can also be described as graded modules over the
  quiver algebra $R = R_{\Gamma,\cA}$ generated by idempotents
  $e_\alpha$, $\alpha\in \Gamma$ in degree $0$, maps
  $\psi_{\beta,\alpha}$, $\beta \ge \alpha$ in degree $1$, together
  with all the elements of the rings $\cA_\alpha$, $\alpha\in \Gamma$,
  and satisfying obvious relations (for instance, for $a\in
  \cA_\gamma$, $\psi_{\beta,\alpha}a = \phi_{\beta,\alpha}(a)$ if
  $\gamma = \alpha$, and is zero otherwise).  $\cA\mof$ is then the
  category of finitely generated $R$-modules.
\end{rmk}

%Again
%let $D(\cA \text{-Mod})$ and $D^b(\cA \text{-Mod})$ denote the corresponding derived
%categories.

\subsubsection{Co-sheaves of modules on a ringed quiver}

Given a ringed quiver $(\Gamma ,\cA _{\Gamma})$, by a {\it co-sheaf}
of $\cA _{\Gamma}$-modules we mean the following data: for every
$\alpha \in \Gamma$ there is given a $\cA _{\alpha}$-module
$\cM_{\alpha}$ with a morphism $\varphi _{\alpha \beta}\colon
\cM_{\alpha}\to \cM_{\beta}$ of $\cA _{\beta}$ modules if $\beta \geq
\alpha$, so that $\varphi _{\beta \gamma}\varphi _{\alpha
  \beta}=\varphi _{\alpha \gamma}$ for $\gamma \geq \beta \geq
\alpha$.  We call co-sheaves of $\cA $-modules co-$\cA $-modules and
denote this abelian category by co-$\cA\Mod$ (resp. co-$\cA\mod$ in
the graded case).

\subsection{DG ringed quiver} \label{DG sheaves}
A {\it sheaf of DG algebras} $\cC=\cC _{\Gamma}$
on $\Gamma$ is defined in the same way as a sheaf of rings, except the
stalks $\cC _{\alpha}$ are DG algebras and morphisms $\phi _{\beta
  \alpha}\colon \cC _{\beta}\to \cC _{\alpha}$ are homomorphisms of DG
algebras. Similarly, a DG $\cC$-module $\cN$ is a collection $\{
\cN_{\alpha}, \psi _{\beta \alpha}\} _{\alpha \leq \beta }$, where
$\cN_{\alpha}$ is a DG module over $\cC _{\alpha}$ and $\psi _{\beta
  \alpha}\colon \cN_{\beta}\to \cN_{\alpha}$ is a homomorphism of DG
modules over $\cC _{\beta}$. Denote by $\cC \text{-DG-Mod}$ the
abelian category of DG $\cC$-modules. One can define a natural
triangulated category $D(\text{DG-$\cC$})$ which is called the derived
category of DG $\cC$-modules (see \cite{L}).

As mentioned before it is sometimes convenient to consider $\Gamma $
as a topological space.  Then $\cC$ induces a sheaf of DG algebras on
this space and DG $\cC$-modules become sheaves of DG modules over that
sheaf of DG algebras.

\begin{defn} Let $\cM$ be a sheaf (or a co-sheaf, or a DG-module) on a quiver $\Gamma$.
  
  a) If $\Phi\subset \Gamma$ is a locally closed subset (i.e.\ the
  difference of two open sets), denote by $\cM _{\Phi}$ the extension
  by zero to $\Gamma$ of the restriction $\cM \vert _{\Phi}$.  Thus
  $(\cM _{\Phi})_\alpha = \cM_\alpha$ if $\alpha \in \Phi$, and $0$
  otherwise.  The restriction map $(\cM _{\Phi})_\alpha \to (\cM
  _{\Phi})_\beta$ is the one from $\cM$ if $\alpha$, $\beta\in \Phi$,
  and is zero otherwise.
  
  b) In case $\cM$ is graded and $k\in \Z$ denote by $\cM\lb k\rb$ 
  the same object
  shifted ``down" by $k$, i.e.  $\cM\lb k\rb_i=\cM _{k+1}$.
\end{defn}

\begin{ex} \label{forget the grading}
  Let $(\Gamma ,\cA )$ be a graded ringed quiver.  Assume that the
  algebras $\cA_\alpha$ are {\em evenly} graded.  We may consider
  $\cA$ as a sheaf of DG algebras with zero differential. There is a
  natural ``forgetful" exact functor between the corresponding derived
  categories
  $$\nu \colon D(\cA\mod)\to D(\text{DG-$\cA$}).$$
  Namely, an object
  of $D(\cA\mod)$ is a complex $\cM^\udot$ of graded $\cA$-modules
  (thus it has a double grading), whereas an object of
  $D(\text{DG-$\cA$})$ is a single DG $\cA$-module. We put $\nu
  (\cM^\udot)=\oplus _i\cM^i\lb -i\rb$ as $\cA$-modules, 
  with the obvious induced differential.  It is easy to see
  that $\nu$ is a triangulated grading (\S\ref{triangulated gradings}).
%Hence the triangulated category $D(\cA\mof)$ may be considered as a {\it mixed}
%version of the triangulated category  $D(\text{DG-$\cA$})$.
\end{ex}

\subsection{Generalities on tori and toric varieties}

Fix a complex torus $T\simeq (\bbC ^*)^n$ with Lie algebra
$\mathfrak{t}$. 
The lattice $N=N_T = \Hom(\bbC^*, T)$ of
co-characters embeds naturally into $\mathfrak{t}$. 
Namely, given a group homomorphism $\phi \colon \bbC^*\to T$, 
the corresponding point in $\mathfrak{t}$ is $d\phi
_*(1)$, where $d\phi _* \colon \bbC \to \mathfrak{t}$ is the induced
map of Lie algebras.  Then $\mathfrak{t} \cong N_\C := N \otimes_Z \C$.

Note that $N$ is naturally isomorphic to the fundamental
group $\pi _1(T)$: given an element $n\in N$ the corresponding map
$f\colon [0,1]\to T$ is defined by the formula
$$f(t)=e^{2\pi i tn}.$$
Clearly this correspondence is functorial with
respect to homomorphisms of tori.

There is also a natural lattice of characters
$M=M_T\subset \mathfrak{t}^*$ defined similarly.  The abelian groups
$M$ and $N$ are dual to each other: $M = \Hom(N,\Z)$. 
The dual torus $\check T$ is the torus for which
$M_{\check T}=N_T$ and $N_{\check T}=M_T$; it 
is isomorphic to $T$, but not canonically.

\subsection{Toric varieties and fans}
Let $X$ be a normal $T$-toric variety.  The $T$-orbits $\{O_{\alpha}\}$
in $X$ are 
indexed by the cones $\alpha$ in a finite polyhedral fan $\Sig$ in
the vector space $N_\R := N \otimes_\Z \R$ which is rational with respect to the
lattice $N$.  If $\alpha \in \Sig$, then 
$\Span_\C(\alpha) \subset \mathfrak{t}$ is 
the Lie algebra of the stabilizer $T_\alpha$ of the corresponding 
orbit $O_\alpha$ (since $T$ is abelian, the stabilizer can be taken
at any point).

We put the 
natural inclusion order on $\Sig$, where
 $\alpha \le \beta$ if and only if $\alpha$ is a face of $\beta$.
Then $\alpha \le \beta$ if and only if $O_\beta \subset \ol{O_\alpha}$.
Thus open unions of orbits correspond to subfans of $\Sig$. 

More generally we will want to consider locally closed unions of
orbits in $X$.  Such a subvariety $Y$ corresponds to a locally closed
subset $\Lambda \subset \Sig$, which satisfies 
all the fan properties, except that it need not be closed under 
taking faces.  Instead if $\alpha \le \beta$ are cones in 
$\Lambda$, then $\Lambda$ must contain all faces of $\beta$ which 
contain $\alpha$.  We call such subsets {\em quasifans}.

\subsection{} \label{toric projections}
Let $X$ be a $T$-toric variety.   
For each orbit $O_\alpha$ in $X$ denote by
$$\St(O_\alpha)=\bigcup_{O_\alpha \subset \overline{O_\beta}}O_\beta$$
its star in $X$.  Consider the stabilizer $T_\alpha \subset T$ of the
orbit $O_\alpha$. Since $X$ is normal, the group $T_\alpha $ is
connected, and hence is a torus.  There exists a non-canonical
homeomorphism
$$\St(O_\alpha)\simeq X^\prime \times (T/T_\alpha) \simeq X' \times
O_\alpha,$$
where $X^\prime $ is an affine $T_\alpha$-toric variety
with a single fixed point.

The action of $T_\alpha $ on $\St (O_\alpha)$ defines a canonical
projection
$$p_\alpha\colon \St(O_\alpha)\to O_\alpha,$$
which is compatible with
the product decomposition above.

If $O_\beta \subset \St(O_\alpha)$ we denote by $p_{\beta \alpha}$ the
restriction of $p_{\alpha}$ to $O_\beta$. The collection of
projections $\{p_{\beta \alpha}\}$ is compatible in the sense that
$p_{\beta \alpha}p_{\gamma \beta}=p_{\gamma \alpha}$ wherever all maps
are defined.

\begin{lemma} \label{fund nbds}
  Each point in $O_\alpha$ has a fundamental system of distinguished
  contractible neighborhoods $U\subset X$, such that $U\cap
  \St(O_\alpha)\subset p_\alpha ^{-1}(U\cap O_\alpha)$ and for each
  orbit $O_\beta \subset \St(O_\alpha)$ the inclusion $U\cap
  O_\beta\hookrightarrow p^{-1}_{\alpha \beta}(U\cap O_\alpha)$ is a
  homotopy equivalence.
\end{lemma}

\subsection{The ringed quiver $(\Sig ,\cA _\Sig)$} 
\label{the sheaf cA}
Let $\Sig$ be a finite polyhedral fan, rational or not, with the 
inclusion partial order on faces.
There is a natural graded ringed quiver $\cA =\cA
_\Sig$ on $\Sig$: for $\tau \in \Sig$ the stalk $\cA
_{\tau}$ is the graded ring of complex-valued polynomial functions on
the span of $\tau$.  The structure homomorphisms $\phi $ are the
restrictions of functions. We consider linear functions as having
degree $2$, so that $\cA $ is evenly graded.

Note that the ringed quiver $\cA$ makes sense even for non-rational
fans (which do not correspond to toric varieties).

\begin{rmk} Notice that the topological space associated to the 
  partially ordered set $\Sig$ is homeomorphic to the quotient space
  $\overline{X}=X/T$.  If $Y\subset X$ is $T$-invariant and locally
  closed, the space of sections of $\cA$ on $\ol{Y}$ is canonically
  identified with the equivariant cohomology $H^*_T(Y;\C)$.
  This is why this sheaf is useful for studying the equivariant
  topology of $X$.  Later in \S\ref{equivariant sheaves} we will use the 
  categories $\cA_X\mof$, DG-$\cA_X$, and their derived categories to 
  model equivariant sheaves and complexes on $X$.
\end{rmk}

\subsection{The ringed quiver $(\Sig ^\circ,\cB _\Sig)$}
Now consider the partially ordered set $\Sig ^\circ$ which is 
$\Sig$ with the opposite ordering.  One may think about $\Sig
^\circ$ as the partially ordered set of orbits of $X$, where $O_\alpha
\leq O_\beta$ iff $O_\alpha\subset \overline{O_\beta}$.  There is a
natural sheaf of rings $\cB _X=\cB$ on $\Sig ^\circ$: for an orbit
$O_\alpha$ take $\cB _{\alpha}$ to be the group ring $\bbC [\pi
_1(O_\alpha)]$. If $O_\alpha \leq O_\beta$ the canonical projections
$p_{\alpha \beta}\colon O_\beta \to O_\alpha$ induce homomorphisms
$\cB _\beta \to \cB _\alpha$.  Thus we obtain a ringed quiver $(\Sig
^\circ,\cB)$. 

As with the ringed quiver $\cA$, $\cB$ can be described entirely 
in terms of the fan $\Sig$, without reference to the 
toric variety.  
For any orbit $O_\alpha$, there is a canonical identification
$\pi_1(O_\alpha) = N_\alpha$, where $N_\alpha$ is the lattice  
$N/(N\cap \Span(\alpha))$.  
If $\alpha$ is a face of $\beta$, so $\beta \le \alpha$
in $\Sig^\circ$, the homomorphism $\pi_1(O_\alpha) \to \pi_1(O_\beta)$
comes from the natural map $N_\alpha \to N_\beta$.

\subsection{The ringed quiver $(\Sig ^\circ, \cT _\Sig)$}
We can define another ringed quiver on $\Sig^\circ$ as
follows. For $\alpha\in \Sig$, let $\cT_\alpha=\Symm(N_{\alpha,\C})$. 
We consider it as a graded
polynomial algebra, where elements of $N_{\alpha,\C}$ have
degree $2$.  If  $\alpha$ is a face of $\beta$, 
the homomorphism $\cT_\alpha \to \cT_\beta$
comes from the natural map $N_\alpha \to N_\beta$.  

Notice that $N_{\alpha,\C}$
is canonically isomorphic to the Lie algebra of $T/T_\alpha$,
which can be canonically identified with $O_\alpha$. 
Thus our graded ringed quiver can be described 
more geometrically:  if $O_\alpha \leq O_\beta$, 
the canonical projections
$p_{\alpha \beta}\colon O_\beta \to O_\alpha$ induce morphisms of tori
$T/T_\beta \to T/T_\alpha$, hence they induce homomorphisms of graded
polynomial algebras $\cT _\beta \to \cT _\alpha$. 

\subsection{Dual affine toric varieties} \label{dual toric varieties}
Let $X$ be an affine $T$-toric variety with a single fixed point.  The
corresponding fan $\Sig = \Sig _X$ consists of a single full-dimensional cone 
$\sigma = \sigma_X\subset N_{T,\bbR}$ together with its faces.
We have the dual cone $\check\sigma$ in the dual vector
space $M_{T,\bbR} = (N_{T,\bbR})^*$, defined by
\[\check\sigma = \{y \in M_{T,\bbR} \mid \la x, y\ra \ge 0\;\text{for all}\;x\in \sig.\]

Let $\check T$ be the dual torus to $T$; then 
$N_{\check T, \bbR} = (N_{T,\bbR})^*$ canonically.  

\begin{defn}
  The dual toric variety $\check X$ to $X$ is the 
  affine $\check T$-toric variety defined by the fan
  $\Sig^\vee$ consisting of $\check \sig$ and all its
  faces, with respect to the lattice $M_T = N_{\check T}$.
  In other words, we have $\check \sigma _X=\sigma_{\check X}$.
\end{defn}

There is an order-reversing isomorphism $\alpha \mapsto \alpha^\bot$
between $\Sig$ and $\Sig^\vee$, defined by $\alpha^\bot = 
\sig^\vee \cap \Span(\alpha)^\bot$.  In particular we have
$\Span(\alpha^\bot) = \Span(\alpha)^\bot$. 
This map gives
an identification $\Sig^\vee = \Sig^\circ$ of partially ordered sets.

With respect to this identification, the 
ringed quivers $(\Sig,\cA_\Sig)$ and 
$((\check\Sig)^\circ,\cT _{\check\Sig})$ are identical.
This will be important 
for the definition of our equivariant-constructible 
duality in \S\ref{main proofs}.

\section{Locally constant sheaves on toric varieties}

\subsection{Some lemmas about sheaves on toric varieties}
For a topological space $Y$ denote by $\Sh(Y)$ the abelian category of
sheaves of complex vector spaces on $Y$.

Let $X$ be a normal toric variety.  We consider $X$ as a topological
space in the classical topology.  Let $Z$ be a $T$-invariant subspace
of $X$. Denote by $LC(Z)\subset \Sh(Z)$ the full subcategory of
sheaves which are locally constant on each orbit.

Consider the full subcategory $D^b_{LC}(\Sh(Z))$ of the bounded
derived category $D^b(\Sh(Z))$, consisting of complexes with
cohomologies in $LC(Z)$.

Fix an orbit $O_\alpha \subset X$ and choose a locally closed
$T$-invariant subset $W\subset \St(O_\alpha)$, which contains
$O_\alpha$. Denote by $i\colon O_\alpha \hookrightarrow W$ the
corresponding closed embedding. Let $j\colon U\hookrightarrow W$ be
the complementary open embedding of $U=W-O_\alpha$. Denote by $q\colon
W\to O_\alpha$, $p\colon U\to O_\alpha$ the restrictions of the
projection $p_\alpha$ to $W$ and $U$ respectively.

\begin{lemma} \label{a lemma}
  In the above notation the functors ${\bf R} q_*$ and $i^*$ from
  $D^b_{LC}(W)$ to $D^b_{LC}(O_\alpha)$ are naturally isomorphic.  In
  particular, the functors ${\bf R} p_*$ and $i^*{\bf R} j_*$ from
  $D^b_{LC}(U)$ to $D^b_{LC}(O_\alpha)$ are naturally isomorphic.
  Hence, the functors $p_*$ and $i^*j_*$ from $LC(U)$ to
  $LC(O_\alpha)$ are naturally isomorphic.
\end{lemma}

\begin{proof} Using distinguished neighborhoods  of points in
  $O_\alpha$ (Lemma \ref{fund nbds}) we see that there exists a
  natural morphism of functors ${\bf R} q_*\to i^*$. Let us show that
  it is an isomorphism.
  
  The category $D^b_{LC}(W)$ is the triangulated envelope of objects
  ${\bf R} j_{\beta *}L$, where $j_\beta \colon O_\beta
  \hookrightarrow W$ is the embedding of an orbit $O_\beta$ and $L$ is
  an object in $LC(O_\beta)$.  So we may assume that $U=O_\beta$ and
  it suffices to show that $i^*{\bf R} j_*L={\bf R}p_*L$ (the case
  $\alpha =\beta$ is clear).
  
  Choose a distinguished neighborhood $V\subset X$ of a point in
  $O_\alpha$. Then the complex $\Gamma (V\cap O_\alpha, i^*{\bf R}
  j_*L)$ is quasi-isomorphic to the complex ${\bf R} \Gamma (V\cap
  O_\beta,L)$. But the inclusion $V\cap O_\beta \subset p^{-1}(V\cap
  O_\alpha)$ is a homotopy equivalence, hence it induces a
  quasi-isomorphism ${\bf R} \Gamma (p^{-1}(V\cap O_\alpha),L)\simeq
  {\bf R} \Gamma (V\cap O_\beta,L)$. This proves that $i^*{\bf R}
  j_*L={\bf R}p_*L$.

  The last statement now follows by taking $H^0$.
\end{proof}

\subsection{Equivalence of derived categories}

\begin{thm} \label{D^b(LC)}
  The natural functor $D^b(LC(X))\to D^b_{LC}(\Sh(X))$ is an
  equivalence.
\end{thm}

\begin{proof} Let $i=i_O\colon O\hookrightarrow X$ be the 
  (locally closed) embedding of an orbit $O$.  For $F\in LC(O)$ we may
  consider two different (derived) direct images of $F$ under the
  embedding $i$: one in the category $D^b_{LC}(\Sh(X))$, denoted as
  usual by ${\bf R}i_*F$, and the other in the category $D^b(LC(X))$,
  which we denote by ${\bf R}_{LC}i_*F$. It is clear that the category
  $D^b_{LC}(\Sh(X))$ (resp. $D^b(LC(X))$) is the triangulated envelope
  of the objects ${\bf R}i_*F$ (resp. ${\bf R}_{LC}i_*F$) for various
  orbits $O$ and locally constant sheaves $F$ on them. So it suffices
  to prove the following two claims.
  
  \medskip

\noindent{\it Claim 1.} The complexes  ${\bf R}i_*F$ and  ${\bf R}_{LC}i_*F$ are
quasi-isomorphic.

\medskip

\noindent{\it Claim 2.} Let $i$ and $F$ be as above, $j\colon O^\prime \hookrightarrow X$
be the embedding of an orbit and $G\in LC(O^\prime)$. Then
$$\Ext ^\udot_{D^b(LC(X))}({\bf R}_{LC}j_*G,{\bf R}_{LC}i_*F)= \Ext
^\udot_{D^b_{LC}(\Sh(X))}({\bf R}_{LC}j_*G,{\bf R}_{LC}i_*F).$$

\medskip

Let us prove the second claim first. Using the adjunction we need to
prove that
$$\Ext ^\udot_{D^b(LC(O))}(i^*{\bf R}_{LC}j_*G,F)= \Ext
^\udot_{D^b_{LC}(\Sh(O))}(i^*{\bf R}_{LC}j_*G,F).$$
By devissage this
is a consequence of the following lemma.

\begin{lemma} Let $Y$ be a $K(\pi ,1)$-space, $LC(Y)$ -- the category of locally constant
  sheaves on $Y$. Then for any $A,B\in LC(Y)$
  $$\Ext ^\udot_{LC(Y)}(A,B)=\Ext ^\udot _{\Sh(Y)}(A,B).$$
\end{lemma}

\begin{proof} Let $f\colon \tilde{Y}\to Y$ be the universal covering map. Then the functor $f^*$
  establishes an equivalence of abelian categories
  $$f^*\colon \Sh(Y)\rightarrow \Sh _{\pi}(\tilde{Y}),$$
  where
  $\Sh_{\pi}(\tilde{Y})$ is the category of $\pi$-equivariant sheaves
  on $\tilde{Y}$ \cite{Gr}. Clearly $f^*$ preserves locally constant
  sheaves.

\begin{rmk}
  It is well known that the category $\Sh_{\pi ,LC}(\tilde{Y})$ of
  locally constant (=constant) $\pi $-equivariant sheaves on
  $\tilde{Y}$ is equivalent to the category of $\pi$-modules. The
  equivalence is provided by the functor of global sections $\Gamma$.
\end{rmk}

Put $\tilde{A}=f^*A$, $\tilde{B}=f^*B$.  We will prove that
$\Ext^\udot_{\Sh_{\pi, LC}(\tilde{Y})}(\tilde{A},\tilde{B})= \Ext
^\udot_{\Sh_{\pi}(\tilde{Y})}(\tilde{A},\tilde{B}).$ Choose an
injective resolution
$$\tilde{B}\to I^0\to I^1\to ...$$
in the category
$\Sh_{\pi,LC}(\tilde{Y})$. It suffices to prove that
$\Ext^k_{\Sh_{\pi}(\tilde{Y})}(\tilde{A},I^t)=0$ for any $t$ and
$k>0$.  Put $I=I^t$ and choose a resolution $0\to I\to J^0\to
J^1\to...$, where $J$'s are injective objects in
$\Sh_{\pi}(\tilde{Y})$. So
$$\Ext ^k_{\Sh_{\pi}(\tilde{Y})}(\tilde{A},I)= H^k(\Hom
^\udot_{\Sh_{\pi}(\tilde{Y})}(\tilde{A},J^\udot)).$$

Notice that $I$, as a sheaf, is constant and each $J^s$, as a sheaf,
is injective \cite{Gr}. Hence the complex of global sections
$$0\to \Gamma (I)\to \Gamma (J^0)\to \Gamma (J^1)\to ...$$
is exact
($\tilde{Y}$ is contractible). Since $\tilde{A}$, as a sheaf, is also
constant, the complex
$$0\to \Hom _{\Sh_{\pi}(\tilde{Y})}(\tilde{A},I)\to
\Hom_{\Sh_{\pi}(\tilde{Y})}(\tilde{A},J^0)\to...$$
is isomorphic to the complex
$$0\to \Hom _{\pi}(\Gamma (\tilde{A}),\Gamma (I))\to \Hom
_{\pi}(\Gamma (\tilde{A}),\Gamma (J^0))\to ...$$
The $\pi$-module
$\Gamma (I)$ is injective and so are the $\pi$-modules $\Gamma (J^s)$
for all $s\geq 0$. Hence the last complex is exact. This proves the
lemma and Claim 2.
\end{proof}

Let us prove Claim 1. Choose a locally constant sheaf $I$ on $O$ which
is injective in the category $LC(O)$. It suffices to prove that the
complex of sheaves ${\bf R}i_*I$ is acyclic except in degree zero.
Choose an orbit $O_\alpha \subset \overline{O}$ and let $p\colon O\to
O_\alpha$ be the canonical projection. Fix a distinguished
neighborhood $U\subset X$ of a point in $O_\alpha$ (Remark \ref{fund
  nbds}). By Lemma \ref{a lemma} above the complex $\Gamma (U\cap
O_\alpha,{\bf R}i_*I)$ is quasi-isomorphic to the complex
$${\bf R}\Gamma ( p^{-1}(U\cap O_\alpha), I)= {\bf R}\Hom
^{\udot}(\bbC _{p^{-1}(U\cap O_\alpha)}, I\vert _{p^{-1}(U\cap
  O_\alpha)}).$$
Since the space $p^{-1}(U\cap O_\alpha)$ is $K(\pi
,1)$, and the restriction of the local system $I$ to $p^{-1}(U\cap
O_\alpha)$ remains injective, it follows from the above lemma that the
complex ${\bf R}\Hom ^{\udot}(\bbC _{p^{-1}(U\cap O_\alpha)}, I\vert
_{p^{-1}(U\cap O_\alpha)})$ is acyclic in positive degrees.  This
proves Claim 1 and the theorem.
\end{proof}

\begin{rmk} The key property of toric varieties which is used in the proof of the above
  theorem is that the star of an orbit is homotopy equivalent to the
  orbit itself. For example, the analogue of the above theorem does
  not hold for $\bbP ^1$ which is stratified by two cells: $\bbC $ and
  a point.
\end{rmk}

The category $LC(X)$ has enough injectives: injective objects are sums
of objects of the form $i_*I$, where $i\colon O \hookrightarrow X$ is
an embedding of an orbit and $I\in LC(O)$ is an injective local
system.  Furthermore, $LC(X)$ has finite cohomological dimension, so
objects in $D^b(LC(X))$ can be represented by bounded complexes of
injectives.

Thus if $j\colon Y\hookrightarrow X$ is an embedding of a locally
closed $T$-invariant subspace, we can take derived functors of $j_*$
and $\Gamma_Y$ (sections with support in $Y$), giving functors 
${\bR}_{LC}j_*\colon D^b(LC(Y)) \to D^b(LC(X))$ and 
${\bR}_{LC}\Gamma_Y\colon D^b(LC(X)) \to D^b(LC(X))$.  Define 
$j^!_{LC} = j^*{\bR}_{LC}\Gamma_Y$.

On the other hand, the usual derived functors restrict to functors
${\bR}j_*\colon D^b_{LC}(Y) \to D^b_{LC}(X)$ and 
${\bR}\Gamma_Y\colon D^b_{LC}(X) \to D^b_{LC}(X)$, and we have 
$j^! = j^*\bR\Gamma_Y$.  The following corollary to the proof of Theorem
\ref{D^b(LC)} will be used later when we discuss the intersection
cohomology sheaves.

\begin{cor} \label{R vs R_LC}
  The functors ${\bR}_{LC}j_*$ and $\bR j_*$ are isomorphic under
  the equivalence of Theorem \ref{D^b(LC)}, as are $j^!_{LC}$ and
  $j^!$.
%Let $j\colon Y\hookrightarrow X$ be the embedding of a $T$-invariant subspace.
%Then the functors ${\bf R}j_*$ and ${\bf R}_{LC}j_*$ from $D^b(LC(U))$ to $D^b_{LC}(X)$ are
%isomorphic.
\end{cor}

\begin{proof} 
  Since $D^b(LC(X))$ is generated by injective objects of $LC(X)$, for
  the first claim it will be enough to consider $i_*I$, where $i\colon
  O \hookrightarrow Y$ is the inclusion of an orbit and $I$ is
  injective in $LC(O)$.  Then $\bR_{LC}j_*(i_*I) = j_*i_*I = (j\circ
  i)_*I\simeq R_{LC}(j\circ i)_*I$, since $i_*I$ is injective.  On the
  other hand, we have
\[Rj_*(i_*I) = Rj_*R_{LC}i_*I \simeq Rj_*Ri_*I \simeq R(j\circ i)_*I ,\] 
using Claim 1 from the proof of Theorem \ref{D^b(LC)}.  Applying Claim
1 once more gives $\bR_{LC}j_*(i_*I) \simeq Rj_*(i_*I)$.
  
For the second part, let $i\colon O\hookrightarrow X$ be the inclusion
of an orbit, and let $I\in LC(O)$ be injective; we will show that
$j^!(i_*I)$ and $j^!_{LC}(i_*I)$ are quasi-isomorphic.  The inclusion
$j$ can be factored as the composition of an open embedding and a
closed embedding.  The required isomorphism is obvious when $Y$ is
open, so we can assume that $Y$ is closed.  If $Y$ contains
$\overline{O}$, then $j^!i_*I = j^*i_*I = j^!_{LC}i_*I$.  Otherwise,
we have $j^!_{LC}(i_*I) = j^*\Gamma_Y(i_*I) = 0$, since $i_*I$ is
injective and all nonzero sections of $i_*I$ must contain points of
$O$ in their support.  On the other hand, $i_*I \simeq \bR i_*I$, so
$j^!(i_*I) = 0$ as well.  This completes the proof.
\end{proof}

\subsection{Quiver description of the category $LC(X)$}

Recall the ringed quiver $(\Sig ^\circ,\cB)$ associated with the toric
variety $X$.  We are going to define a functor
$$\eta \colon LC(X)\to \co\cB\Mod.$$
For this we need to recall how to
glue sheaves on topological spaces. Surely this construction is well
known, but we do not know a reference.

Let $Y$ be a topological space, $i\colon Z\hookrightarrow Y$ the
embedding of a closed subset and $j\colon U=Y-Z\hookrightarrow Y$ the
complementary open embedding. Consider the abelian category $\Sh(Y,Z)$
consisting of triples $(G,H,\xi)$, where $G\in \Sh(Z)$, $H\in \Sh(U)$
and $\xi $ is a morphism of sheaves in $\Sh(Z)$ $\xi\colon G\to
i^*j_*H$.  We have a natural functor $\tau \colon \Sh(Y)\to \Sh(Y,Z)$
which associates to a sheaf $F\in \Sh(Y)$ its restrictions $i^*F\in
\Sh(Z)$, $j^*F\in \Sh(U)$ and the pullback under $i^*$ of the
adjunction morphism $F\to j_*j^*F$.

\begin{lemma} The functor $\tau$ is an equivalence.
\end{lemma}

\begin{proof} Let us define the inverse functor $\eta \colon \Sh(Y,Z)\to
  \Sh(Y)$. Given $(G,H,\xi)\in \Sh(Y,Z)$ define a presheaf
  $\overline{F}$ on $Y$ as follows.  For an open subset $V\subset Y$
  put $\overline{F}(V)=H(V)$ if $V\subset U$. Otherwise set
  $$\overline{F}(V)=\{(g,h)\in G(V\cap Z)\times H(V\cap U)\mid \xi
  (g)=h^\prime \},$$
  where $h^\prime $ is the image of $h$ in
  $i^*j_*H(V\cap Z)$. Then let $\eta (G,H,\xi)\in \Sh(Y)$ be the
  sheafification of $\overline{F}$.
\end{proof}

Let $F\in LC(X)$. Denote the stalk of $F$ at the distinguished point
in an orbit $O_\alpha$ by $F_\alpha$. Then $F_\alpha $ is a $\cB
_\alpha $-module. Given two orbits $O_\alpha \subset
\overline{O_\beta}$ consider the canonical projection $p = p_{\beta
  \alpha}\colon O_\beta \to O_\alpha$. By Lemma 6.7 and Lemma 6.1 the
restriction of the sheaf $F$ to the union of the two orbits defines a
morphism of sheaves
$$F\vert _{O_\alpha}\to p_*(F\vert _{O_\beta}),$$
or equivalently a
morphism of sheaves
$$p^{-1}(F\vert _{O_\alpha})\to (F\vert _{O_\beta}).$$
Such a morphism
is equivalent to a homomorphism of $\cB _\beta$-modules $F_\alpha\to
F_\beta$.  So the sheaf $F$ defines a co-$\cB$-module. This is our
functor
$$\eta \colon LC(X)\to \co\cB\Mod.$$

\begin{thm} The functor $\eta$ is an equivalence.
\end{thm}

For example, in case the toric variety $X$ is the affine line with two
orbits, $\bbC ^*$ and the origin, an object in $LC(X)$ is the same as
a vector space $P$, a $\bbC [\pi _1 (\bbC ^*)]$-module $Q$ and a linear 
map $P\to Q^{\pi_1(\bbC^*)}$.

\begin{proof} We will prove the theorem by induction on the number of
  orbits in $X$. For one orbit the statement of the theorem is a well
  known equivalence between the category of locally constant sheaves
  and that of $\pi _1$-modules.
  
  Now we proceed with the induction step. Pick an orbit $O_\alpha
  \subset X$ of smallest dimension. We may assume that
  $X=\St(O_\alpha)$. Indeed, otherwise $X$ may be covered by open
  $T$-invariant subsets $V$, which are strictly smaller than $X$. By
  induction, the theorem is true for each $V$ and so we obtain the
  equivalence for $X$ by gluing the corresponding equivalences
  $LC(V)\simeq \co\cB _V\Mod$.
  
  Put $U=X-O_\alpha$ and let $j\colon U\hookrightarrow X$ and $i\colon
  O_\alpha \hookrightarrow X$ be the open and closed embeddings
  respectively.  By Lemma 6.7 a sheaf $F\in LC(X)$ is the same as a
  triple $(G,H,\xi)$, where $G\in LC(O_\alpha)$, $H\in LC(U)$ and $\xi
  \colon G\to i^*j_*H$. By Lemma 6.1 $i^*j_*H=p_{\alpha *}H$. Thus by
  adjunction the morphism $\xi$ is the same as a morphism $\iota
  \colon p_\alpha^{-1}G\to H$. Let $G_\alpha$ be the $\cB
  _\alpha$-module corresponding to $G$. It is easy to see that the
  sheaf $p_\alpha^{-1}G$ considered as a co-$\cB _U$-module is the
  constant one equal to $G_\alpha$. Thus the triple $(G,H,\iota)$ is
  the same as a co-$\cB$-module. This proves the theorem.
\end{proof}

\subsection{Unipotent sheaves}

\begin{defn} A sheaf $F\in LC(X)$ is called {\it unipotent} if for each
  orbit $O_\alpha$ and $x\in \pi_1(O_\alpha)$ the action of the
  operator $x-1$ on the stalk $F_\alpha$ of $F$ at a point of
  $O_\alpha$ is locally nilpotent. It is called {\it co-finite} if in
  addition the space of invariants $F_\alpha ^{\pi _1(O_\alpha)}$ is
  finite-dimensional for all $\alpha$. Let $LC_{u}(X)$ and
  $LC_{cf}(X)$ be the full subcategories of $LC(X)$ consisting of
  unipotent (resp.\ co-finite) sheaves.  

  The next result describes the corresponding subcategories of 
  $\co\cT\Mod$ under the equivalence $\eta$.  Let 
  $\co\cT\Mod_{n}$ be the full subcategory of co-$\cT$-modules 
  $\cM$ which are ``supported at the origin'', i.e.\ for which every
$m\in \cM_\alpha$ is annihilated by some power of the homogeneous 
maximal ideal $\m_\alpha\subset \cT _\alpha$.  Let $\co\cT\Mod_{cf}$
be the further full subcategory of modules $\cM$ for which each 
$\cM_\alpha$ is a cofinite $\cT_\alpha$-module (\S\ref{module conventions}).
In other words, in addition to being supported at the origin, 
for each $\alpha$ the space 
$\{m\in \cM_\alpha\mid \m_\alpha \cdot m = 0\}$ should be finite dimensional.
\end{defn}

\begin{thm} \label{cofinite LC sheaves}
The functor $\eta$ restricts to give equivalences of full
abelian subcategories
  $$LC_u(X)\simeq \co\cT\Mod _n,\;\text{and}$$
  $$LC_{cf}(X)\simeq \co\cT\Mcf.$$
\end{thm}

\begin{proof}  Take a sheaf $F \in LC_u(X)$, and let 
$\cM = \eta(F) \in \co\cB\Mod$ be the corresponding $\co\cB$-module.  
Since the action of any
$x\in \pi _1(O_\alpha)\cong N_\alpha$ on $\cM_\alpha$ is
unipotent, the action of the power series
$\frac{1}{2\pi i}\ln x$ is well-defined.  
We can extend this uniquely to a map 
$v \mapsto \frac{1}{2\pi i}\ln v$ 
from $N_{\alpha} \otimes \C$ to $\End(\cM_\alpha)$.  
Any two of these operators commute, since $\pi_1(O_\alpha)$ is abelian.

This gives $\cM_\alpha$ the structure of a $\Sym(N_{\alpha}\otimes\C) =
\cT_\alpha$-module, and in fact makes $\cM$ into a co-$\cT$-module.
The resulting co-$\cT$-module is clearly supported at the origin.
Conversely, given a co-$\cT$-module $\cM$ supported at the origin, 
we can exponentiate the action of elements of $N_\alpha$
to get an action of $\pi_1(O_\alpha)$ on $\cM_\alpha$.  These 
actions combine to give the structure of a $\co\cB$-module
on $\cM$, and then applying $\eta^{-1}$ gives the required object
in $LC_u(X)$.

The second equivalance follows immediately, since the maximal
ideal $\m_\alpha \subset \cT_\alpha$ is generated by 
$\frac{1}{2\pi i}\ln x, x\in \pi_1(O_\alpha)$.
\end{proof}

\section{Mixed locally constant sheaves}
 
\subsection{pre-$\cF$-sheaves}
For toric varieties the Frobenius endomorphism has a natural lift to
characteristic zero -- see \cite{We}.  We will use it to define a 
mixed version of the
category $LC(X)$.

Consider the group homomorphism $\phi \colon T\to T$, $a\mapsto a^2$.
For any toric variety $X$ the homomorphism $\phi$ extends uniquely to
a morphism $\cF=\cF_X\colon X\to X$. Namely, recall that each orbit
$O_\alpha$ is identified with the quotient torus $T/T_\alpha$; then
the map $\cF\colon O_\alpha \to O_\alpha$ is again squaring. The maps
$\phi$ and $\cF$ have degree $2^n$.

\begin{defn} A pre-$\cF$-sheaf is a pair 
$(F,\theta)$, where $F\in LC(X)$ and $\theta $
  is an isomorphism
  $$\theta \colon \cF ^{-1}F\to F .$$
%  We denote the abelian category
%  of pre-$\cF$-sheaves by pre-$LC_\cF(X)$.
\end{defn}

Let us describe the inverse image functor $\cF ^{-1}\colon LC(X)\to
LC(X)$ in terms of co-$\cB$-modules.  The map $\cF$ induces the
endomorphism of the sheaf $\cB$, where each element $x\in
\pi_1(O_\alpha)$ maps to $x^2$. Denote this endomorphism $\psi \colon
\cB \to \cB$.  Fix $F\in LC(X)$ and let $\cM$ be the corresponding
co-$\cB$-module. Then the sheaf $\cF ^{-1}F$ corresponds to the
co-$\cB$-module $\psi _*\cM$, i.e. it is obtained from $\cM$ by
restriction of scalars via $\psi$.  So the isomorphism $\theta \colon
\psi_*\cM\to \cM $ corresponding to the isomorphism $\theta \colon \cF
^{-1}F\to F$ amounts to a compatible system of linear maps $\theta
_\alpha \colon \cM_\alpha \to \cM_\alpha $ such that for $x\in \pi
_1(O_\alpha)$, $m\in \cM_\alpha$
$$\theta _\alpha(x^2m)=x\theta _\alpha (m).$$

\subsection{$\cF$-sheaves and graded co-finite co-$\cT$-modules}
\label{F-sheaves}
\begin{defn} A pre-$\cF$-sheaf $(F,\theta)$
  is called an $\cF$-{\it sheaf} if $F$ is co-finite (unipotent) and
  for each $\alpha$ the endomorphism $\theta _\alpha \colon F_\alpha
  \to F_\alpha $ is diagonalizable with eigenvalues $2^{n/2}$, $n\in
  \bbZ$.  We will refer to $\cF$-sheaves
  on a single orbit as ``$\cF$-local systems''. 
  Denote by $LC_\cF(X)$ the category of $\cF$-sheaves on $X$.
\end{defn}
For any $n\in \Z$ define an automorphism $\la n \ra$ of $LC_\cF(X)$  
by $(F,\theta) \mapsto (F,2^{n/2}\theta)$.

Let $(F,\theta)$ be an $\cF$-sheaf, and take $x\in \pi_1(O_\alpha)$. 
The relation
$$\theta \cdot x^2=x\cdot \theta$$
in $\End(F_\alpha)$ is equivalent
to
$$\theta \cdot 2(\frac{1}{2\pi i}\ln x)=(\frac{1}{2\pi i}\ln x)\cdot
\theta.$$
Thus $F$ considered as the co-$\cT$-module $\cM$ via Theorem
\ref{cofinite LC sheaves} is graded:
$$\cM_k=\{ a\in \cM\mid \theta (a)=2^{-k/2}a\},$$
and the operators
$\frac{1}{2\pi i}\ln x$ map $\cM_k$ to $\cM_{k+2}$.

\begin{defn}
  Let co-$\cT\mod _n$ denote the category of {\it graded}
  co-$\cT$-modules which are supported at the origin, and let
  co-$\cT\Mcf$ denote the full subcategory of objects
  $\cM$ for which $\cM_\alpha$ is a co-finite graded 
  $\cT_\alpha$-module for all $\alpha$.
\end{defn}

The following result is an immediate consequence of Theorem
\ref{cofinite LC sheaves} and the above discussion.

\begin{thm} \label{combinatorial F-sheaves}
  There is a natural equivalence of abelian
  categories
  $$LC_\cF(X)\simeq \co\cT\mcf,$$
  and hence an equivalence
  $$D^b(LC_\cF(X))\simeq D^b(\co\cT\mcf).$$\
Under these equivalences the twist operator $\la n\ra$ goes to 
the grading shift $\lb n\rb$.
\end{thm}

Note that these isomorphisms and the isomorphisms of Theorem
\ref{cofinite LC sheaves} are compatible with the forgetful functors
$LC_\cF(X) \to LC_{cf}(X)$ and co-$\cT\mcf \to $ co-$\cT\Mcf$.  This
means that $D^b(LC_\cF(X)) \to D^b(LC_{cf}(X))$ is a
triangulated grading in the sense of 
\S\ref{triangulated gradings}.

\subsection{Simple and injective mixed sheaves} \label{injective mixed sheaves}
Since the category co-$\cT\mcf$ has enough injectives, so does
$LC_\cF(X)$.  It will be helpful to have a concrete description of
simple and injective objects in this category.

First consider the case of a single $T$-orbit $O = O_\alpha$.  
We have an equivalence 
$LC_\cF(O_\alpha) \simeq \co\cT_\alpha\mcf$.  Up to degree
shifts there is one simple object of $\co\cT_\alpha\mcf$, namely
$(\cT_\alpha/\m_\alpha\cT_\alpha)^*$.  The corresponding object in
$LC_\cF(O_\alpha)$ is the constant local system $\C_{O_\alpha}$, with
$\cF$-structure given by $\theta_\alpha = 1$.  We will denote this
$\cF$-local system by $\C_\alpha$.

The injective hull of $(\cT_\alpha/\m_\alpha\cT_\alpha)^*$ is
$\cT_\alpha^*$; let $\Theta_\alpha$ denote the corresponding injective
object in $LC_\cF(O_\alpha)$.  It has the following topological
description.  Let $q_\alpha\colon \wt{O}_\alpha \to O_\alpha$ be the
universal cover of $O_\alpha$.  Then $\Theta_\alpha$ is the largest
subsheaf of the local system $q_{\alpha
  *}\C_{\wt{O}_\alpha}$ on which all the monodromy operators $x\in
\pi_1(O_\alpha)$ act (locally) unipotently.

Let $b \in O_\alpha$ be the distinguished point.  Since there is a canonical
identification $\pi_1(O_\alpha) \cong N_\alpha$, we can identify the
stalk $(q_{\alpha *}\C_{\wt{O}_\alpha})_b$ with the space of functions
$N_\alpha \to \C$, at the price of choosing a point $\tilde b \in
q^{-1}(b)$.  The action of $x \in \pi_1(O_\alpha)$ on this stalk is 
identified with the pushforward by
the translation $\tau_x\colon n\mapsto n + x$ of the lattice
$N_\alpha$.

The stalk $(\Theta_\alpha)_b$ is thus the space of all functions
$N_\alpha \to \C$ which are annihilated by some power of $x - 1$
for every $x \in N_\alpha$.  This is the space of polynomial
functions $N_\alpha \to \C$.  The logarithm of $x$ acts on these
functions as the differential operator $\partial_x$.

We can now make $\Theta_\alpha$ into an $\cF$-sheaf by letting the
operator $\theta_\alpha$ be the pullback by $N_\alpha \to N_\alpha$,
$x\mapsto 2x$.  The corresponding grading is just our usual even
grading on polynomial functions.  The resulting $\cF$-sheaf is the
injective hull of $\C_\alpha$.

If $X$ has more than one orbit, then up to grading shift the injective
objects of co-$\cT\mcf$ are the sheaves $\cT^*_{\ol \alpha}$ for
$\alpha\in \Sig_X$, where the closure is taken in the fan topology, so
$\ol{\alpha} = \{\beta\in \Sig_X \mid \alpha \le \beta\}$.

The corresponding injective objects in $LC_\cF(X)$ are (up to twists
$\la n\ra$) the sheaves $j_{\alpha *} \Theta_\alpha, \alpha\in \Sig_X$
where $j_\alpha\colon O_\alpha \to X$ is the inclusion.  $j_{\alpha *}
\Theta_\alpha$ is the injective hull of the extension by zero
$j_{\alpha!}\C_\alpha$ of $\C_\alpha$.

Note that the forgetful functor $F_{cf}\colon LC_\cF(X) \to
LC_{cf}(X)$ preserves injectivity, as does the inclusion 
$LC_{cf}(X) \subset LC(X)$.  In particular, this implies that
$D^b(LC_{cf}(X))$ is a full subcategory of $D^b(LC(X))$.  
 
\subsection{Extension and restriction functors}
Let $j\colon Y \hookrightarrow X$ be the inclusion of a $T$-invariant
locally closed subset of $X$.  Since $LC_\cF(X)$, $LC_\cF(Y)$ have
enough injectives, we can take derived functors of the left exact 
functors $j_*$ and $j^*\Gamma_Y$ to get functors
$\bR j_* \colon D^b(LC_\cF(Y)) \to D^b(LC_\cF(X))$ and
$j^!\colon D^b(LC_\cF(X)) \to D^b(LC_\cF(Y))$.  The restriction 
and extension by zero functors $j_!$ and $j^*$ are already exact,
so they do not need to be derived.  

In the same way we get derived functors between 
$D^b(LC(X))$ and $D^b(LC(Y))$.  We will denote them by the same symbols
$\bR j_*$, $j_!$, $j^*$, $j^!$; context will make clear which functor is
meant.  These functors correspond to the ones on $\cF$-sheaves:
$\bR j_* F_{cf} = F_{cf} \bR j_*$, etc.  Furthermore, by Corollary
\ref{R vs R_LC}, these functors agree with the usual topological versions. 

\subsection{Perverse $t$-structure} \label{perverse t-structure}
These functors satisfy the usual adjuntions and distinguished
triangles which allow one to define perverse $t$-structures; see
\cite{GM} or \cite{Br}.

To do this, define $c(\alpha) = \rank N - \dim \alpha = \dim_\C
O_\alpha$ for any $\alpha\in \Sig$, and define full subcategories of
$D = D^b(LC_\cF(X))$ by
\[D^{\le 0}_\cF(X) = \{F^\udot \in D \mid H^i(j_\alpha^* F^\udot) = 0
\;\text{for}\; i > - c(\alpha)\},\]
\[D^{\ge 0}_\cF(X) = \{F^\udot \in D \mid 
H^i({j^!_\alpha} F^\udot) = 0 \;\text{for}\; i < -c(\alpha)\}.\] The
core $P_\cF(X) = D^{\le 0}_\cF(X) \cap D^{\ge 0}_\cF(X)$ is an abelian
category whose objects will be called perverse $\cF$-sheaves.

The same formulas define a perverse $t$-structure $(D^{\le 0}_{cf}(X),
D^{\ge 0}_{cf}(X))$ on $D^b(LC_{cf}(X))$.  The resulting core of
perverse objects will be denoted $P_{cf}(X)$.  The forgetful functor
$F_{cf}\colon D^b(LC_\cF(X))\to D^b(LC_{cf}(X))$ is $t$-exact, so it
restricts to an exact functor $P_\cF(X)\to P_{cf}(X)$.

Simple objects in $P_\cF(X)$ and $P_{cf}(X)$ are obtained as usual by
applying the Deligne-Goresky-MacPherson middle extension
$j_{\alpha!*}$ to a simple local system on an orbit $O_\alpha$,
shifted so as to be perverse.  In particular, \[L^\udot_\alpha :=
j_{\alpha!*} \C_\alpha[c(\alpha)]\la -c(\alpha)\ra\] is simple in
$P_\cF(X)$, and all simple objects are isomorphic to
$L^\udot_\alpha\la n\ra$ for some $\alpha \in \Sig$, $n\in \Z$ (we add
the twist by $-c(\alpha)$ so $L^\udot_\alpha$ will have weight $0$ in
the mixed structure we define below).  Applying the forgetful functor
$F_{cf}$ to $L^\udot_\alpha$ gives the usual intersection cohomology
sheaf $IC^\udot(\ol{O_\alpha};\C)$; these give all the simple objects
of $P_{cf}(X)$.

Note that unlike the usual category of constructible perverse sheaves,
$P_\cF(X)$ is not artinian, since even if $X$ has only one stratum,
objects like $\Theta_\alpha$ have infinite length.  However, Homs are
finite-dimensional in $D^b(LC_\cF(X))$, and hence in $P_\cF(X)$.

\begin{prop} \label{lower star is exact}
Let $i\colon O \hookrightarrow X$ be the inclusion of an 
  orbit.  Then the functor $\bR i_* \colon D^b(LC_\square(O)) \to
  D^b(LC_\square(X))$ is $t$-exact, $\square = \cF,cf$.
\end{prop}

\begin{proof} Suppose that $O=O_\alpha$.  
%  Since $j^!_\beta\bR i_*S^\udot = 0$ for any $S^\udot \in
%  D^b(LC_\cF(O))$ and any $O_\beta \ne O_\alpha$, 
  We have $\bR i_*(D^{\ge 0}_\cF(O)) \subset D^{\ge 0}_\cF(X)$
  automatically, since $j^!_\beta\bR i_* = 0$ for any $\beta \ne
  \alpha$.
  
  Suppose that $S^\udot \in D^{\le 0}(LC_\cF(O))$, and let $M^\udot
  \in D^b(\cT_\alpha\mcf)$ be the corresponding complex of
  $\cT_\alpha$-modules; we have $H^d(M^\udot) = 0$ for $d >
  -c(\alpha)$.  Using the description of injective $\cF$-sheaves from
  \S\ref{injective mixed sheaves} we see that the complex in
  $D^b(\cT_\beta\mcf)$ corresponding to $j^*_\beta\bR i_*S^\udot$ is
  $M_\beta^\udot = \bR\hom_{\cT_\alpha}(\cT_\beta, M^\udot)$.  The
  functor $\bR\hom$ can be defined by deriving either the first or the
  second variable, so the fact that $H^d(M_\beta^\udot) = 0$ \ for $d
  > -c(\beta)$ follows from the fact that $\cT_\beta$ has a resolution
  of length $c(\alpha) - c(\beta)$ by free $\cT_\alpha$-modules.
\end{proof}

Thus for any $\alpha \in \Sig$ we can define an object in $P_\cF(X)$
by
\[\nabla^\udot_\alpha = \bR j_{\alpha *} 
\Theta_\alpha[c(\alpha)]\la -c(\alpha)\ra.\] Note that since
$\Theta_\alpha$ is an injective $\cF$-local system, taking $j_{\alpha
  *}$ instead of $\bR j_{\alpha *}$ defines the same object.  Under
the isomorphism of Theorem \ref{combinatorial F-sheaves},
$\nabla^\udot_\alpha$ corresponds to
$\cT^*_{\overline\alpha}[c(\alpha)]\lb -c(\alpha)\rb$.

This object will be important in the proof of the main properties of
our Koszul duality functor in \S\ref{main proofs} below.

\subsection{Constructible $\cF$-sheaves}
We can also consider the full subcategories $D^b_c(LC_\cF(X)) \subset
D^b(LC_\cF(X))$ and $D^b_c(LC_{cf}(X)) \subset D^b(LC_{cf}(X))$
consisting of complexes $S^\udot$
whose cohomology sheaves have finite
dimensional stalks on each orbit $O_\alpha$.  We call
such objects ``constructible''.  Note that by Theorem \ref{D^b(LC)},
$D^b_c(LC_{cf}(X))$ is equivalent to a full subcategory of the usual
constructible derived category of $D^b_c(X)$: 
namely the category of objects whose cohomology sheaves are orbit-constructible
(and have finite-dimensional stalks), with 
unipotent monodromy on each orbit.

The $t$-structures we defined in the previous section restrict to
$t$-structures on these subcategories, giving abelian cores
$P_{\cF,c}(X) \subset P_\cF(X)$ and $P_{cf,c}(X) \subset P_{cf}(X)$.

\begin{prop} \label{constructible = finite length}
  $P_{\cF,c}(X)$, (resp. $P_{cf,c}(X)$) is the full subcategory of
  objects in $P_\cF(X)$ (resp. $P_{cf}(X)$) consisting of all objects
  of finite length.  In particular, $P_{cf,c}(X)$ is equivalent to the
  full subcategory of the category of constructible perverse sheaves
  on $X$ consisting of objects all of whose simple constituents are of
  the form $IC^\udot(\ol{O_\alpha};\C)$, $\alpha \in \Sig$.
\end{prop}

\begin{rmk} In \cite{Br} a triangulated category $\bD(\Sig)$
  was defined for any fan $\Sig$ to model mixed $T$-constructible
  complexes on the toric variety $X_\Sig$ (in the case $\Sig$ is
  rational).  It can be shown that $D^b_c(LC_\cF(X))$ is equivalent to
  $\bD(\Sig)$; under this equivalence the $t$-structure here is the
  same as the $t$-structure in \cite{Br}.
\end{rmk}

\subsection{Mixed structure and pure $\cF$-sheaves}
In the categories of mixed $l$-adic sheaves or mixed Hodge modules,
simple perverse objects are pure.  We need the following analog of
this fact in our combinatorial setting.  We call an object $S^\udot
\in D^b(LC_\cF(X))$ {\em pure of weight $0$} if for any orbit
$O_\beta$ and any $i \in \Z$, the $\cF$-local systems
$H^i(j_\beta^*S^\udot)$ and $H^i(j_\beta^!S^\udot)$ are direct sums of
finitely many copies of $\C_{\beta}\la i\ra$.  More generally we say
$F^\udot$ is pure of weight $k$ if $S^\udot[-k]$ is pure of weight
$0$.

\begin{prop} \label{unipotent Koszul}
  If $S_1^\udot, S_2^\udot\in P_{\cF,c}(X)$ are pure of weights $r_1$
  and $r_2$, respectively, then $\Hom_{D^b(LC_\cF(X))}(S_1^\udot,
  S_2^\udot[k]) = 0$ unless $r_2 = r_1 - k$.  In particular,
  $\Ext^1_{P_\cF(X)}(S_1^\udot, S_2^\udot) = 0$ unless $r_2 = r_1 -
  1$.
\end{prop}
\begin{proof} There is a spectral sequence with $E_1$ term
\[E_1^{p,q} = \bigoplus_{\dim\alpha = p} \Hom_{D^b(LC_\cF(O_\alpha))}
(j_\alpha^*S_1^\udot,j_\alpha^!S_2^\udot),\] which converges to
$\Hom_{D^b(LC_\cF(X))}(S_1^\udot, S_2^\udot[p+q])$.  Theorem
\ref{purity of simples} implies that $E_1^{p,q} = 0$ unless $p + q =
r_1 - r_2$, which implies the result.
\end{proof}

\begin{thm} \label{purity of simples}\label{mixed F-sheaves}
  The simple perverse sheaf $L^\udot_\alpha$ is pure of weight $0$.
\end{thm}

The proof will be given in \S\ref{appendix}.

\begin{rmk} Purity of IC sheaves enters our main argument twice,
  once via Theorem \ref{purity of simples}, and once in the next
  section, where equivariant IC sheaves are used. In that setting
  the purity follows from a proof of Karu \cite{Ka}, which makes sense 
  even for non-rational fans.  In fact Theorem \ref{purity of simples} 
  can also be stated and proved
  for non-rational fans, without reference to a toric variety.  
  Although the category $LC_\cF(X)$ doesn't make sense, $\co\cT\mcf$
  still does, and one
  can define a functor from the combinatorial equivariant sheaves
  ($\cA$-modules) considered in the next section to $D^b(\co\cT\mcf)$,
  which sends the equivariant IC sheaves to the $L^\udot_\alpha$'s.
  The required purity can then be deduced from Karu's result.
%A similar purity result was proved in
%\cite[Theorem 4.3.1]{Br} for a different combinatorial category 
%modelling mixed sheaves on toric varieties.
\end{rmk}
 
%\subsection{Mixed structure on $P_\cF(X)$} 
We now define a mixed structure on $P_\cF(X)$.  We have already defined 
the twist functor $\la 1\ra$.  What remains is to construct the
filtration $W_\udot$.  

\begin{thm} \label{weight filtration}
There exists a unique
functorial increasing filtration $W_\udot$ on objects of $P_\cF(X)$ satisfying
the following:
\begin{enumerate}

\item[(a)] For any $S^\udot \in P_\cF(X)$ there exists $n\in \Z$ depending on $S^\udot$ 
so that $W_{n}S^\udot = 0$,
\item[(b)] For all $i$ and all $S^\udot\in  P_\cF(X)$, $\Gr_i^WS^\udot = W_iS^\udot/W_{i-1}S^\udot$ 
  is isomorphic to a finite direct sum
  of objects $L^\udot_\alpha\la i\ra$ 
  (thus $\Gr_i^WS^\udot$ is pure of weight $i$), and 
\item[(c)] $(P_\cF(X),W_\udot,\la 1\ra)$ is a
mixed category (\S\ref{mixed categories}).
\end{enumerate}
\end{thm}

For finite length objects, i.e.\ objects in $P_{\cF,c}(X)$, this
follows in a standard way from Theorem \ref{purity of simples} and
Proposition \ref{unipotent Koszul}. 
We give the complete proof in \S\ref{appendix}.

\begin{cor} \label{weights of nabla}
  For any $\alpha \in \Sig$, we have $W_0\nabla^\udot_\alpha \cong
  L^\udot_\alpha$.  In particular, $\nabla^\udot_\alpha$ has weights
  $\ge 0$.
\end{cor}
\begin{proof} If $m$ is the minimum weight in $\nabla^\udot_\alpha$, then 
  $W_m\nabla^\udot_\alpha$ is a semisimple subobject of
  $\nabla^\udot_\alpha$.  But by adjunction $\Hom(L_\beta^\udot\la
  k\ra, \nabla^\udot_\alpha)$ is one-dimensional if $\alpha = \beta$
  and $k = 0$, and vanishes otherwise.
\end{proof}

\section{Equivariant sheaves} \label{equivariant sheaves}

\subsection{} \label{equivariant sheaves and DG modules}

Let us very briefly recall the notion of the bounded, constructible
equivariant derived category $D_{T}^b(X)$ \cite{BL} (note that in
\cite{BL} this category was denoted $D^b_{T,c}(X)$). Let $E$ be a
contructible space with a free $T$-action, and put $X_T=(X\times
E)/T$. Then $E/T=BT$ is the classifying space for $T$ and $X_T\to BT$
is a locally trivial fibration with the fiber $X$. Similarly, a
$T$-invariant subspace $U\subset X$ induces the corresponding subspace
$U_T\subset X_T$. The triangulated category $D_{T}^b(X)$ can be
canonically identified as a full triangulated subcategory of the
bounded derived categories of sheaves on $X_T$. For example, it can be
defined as the triangulated envelope of the collection of all sheaves
$\{ \bbC _{U_T}\}$ ($\bbC _{U_T}$ is the extension by zero to $X_T$ of
the constant sheaf $\bbC$ on $U_T$), where $U\subset X$ is a star of
an orbit. The following theorem is one of the main results in
\cite{L}.

\begin{thm} \label{equivariant-DG equivalence}
  There exists a natural equivalence of triangulated categories
  $$\epsilon\colon D_{T}^b(X)\to D_f(\text{DG-$\cA_X$}).$$
\end{thm}

The equivalence of categories in the above theorem is of ``local
nature" and actually comes from a continuous map of topological
spaces. Namely there exists a natural map
$$\mu \colon X_T\to X/T,$$
and the functor $\epsilon$ is essentially
the derived direct image functor ${\bf R}\mu _*$. In particular, if
$U\subset X$ is the star of an orbit in $X$ and $\tau \in \Sig$ is the
cone corresponding to that orbit, then $\epsilon(\bbC _{U_T})=\cA
_{[\tau]}$. Also $\epsilon $ takes the constant sheaf $\bbC _{X_T}$ on
$X_T$ to the sheaf $\cA$. (In \cite{L} the sheaf $\cA$ is denoted by
$\cH$).

This allows us to define the functor $F_T\colon D^b(\cA) \to D^b_T(X)$
described in the introduction (\S\ref{intro to main results}).  It is
obtained by composing the functor of Example \ref{forget the grading}
with $\epsilon^{-1}$.

\subsection{Combinatorial equivariant complexes}
Let $\Sigma$ be a fan in the vector space $V$, and let $\cA =
\cA_\Sig$ be the sheaf of conewise polynomial functions introduced in
\S\ref{the sheaf cA}.  For the remainder of this section all our
$\cA$-modules will be assumed to be locally finite, so to simplify
notation we put $D^b(\cA\mof) = D^b(\cA)$.

If $\Delta \subset \Sig$ is a subfan or more generally a difference of
subfans, the space of sections of a sheaf $\cM$ on $\Delta$ will be
denoted $\cM(\Delta)$.  If $\Lambda \subset \Delta$ is another subfan,
we put $\cM(\Delta,\Lambda) = \ker(\cM(\Delta) \to \cM(\Lambda))$.

If $\sig \in \Sigma$, recall that $[\sig]$ is the fan of all faces of
$\sig$.  It follows that $\cM({[\sig]})$ is isomorphic to the stalk
$\cM_{\sig}$.  Define $\bdy\sig = [\sig] \setminus \{\sig\}$.  To
simplify notation, we write $\cM(\sig,\bdy\sig) =
\cM([\sig],\bdy\sig)$.

Note that a map of $\cA$-modules $\cM \to \cN$ is determined by the
collection of induced maps $\cM(\sig)\to \cN(\sig)$ over all cones
$\sig\in \Sig$.  A sequence $\cE\to\cM\to \cN$ is exact if and only if
$\cE(\sig) \to \cM(\sig) \to \cN(\sig)$ is exact for all $\sig\in
\Sig$.

%The main theorem of this section is that $D^b(\cA\mof)$ is isomorphic
%to the homotopy category of complexes of ``combinatorially pure'' modules.  
%This simplifies computations considerably, as there is no need to 
%invert quasi-isomorphisms.

\begin{defn} Let $\cM$ be an $\cA$-module.  If
  the stalk $\cM(\sig)$ is a free $\cA_\sig$-module for every $\sig
  \in \Sig$, we say that $\cM$ is {\em locally free}.  If the
  restriction $\cM(\sig) \to \cM({\bdy \sig})$ is surjective for every
  $\sig \in \Sig$, we say that $\cM$ is {\em flabby}.  If both
  conditions hold, we say $\cM$ is {\em combinatorially pure}
  (``pure'' for short).
\end{defn}
Note that flabbiness of $\cM$ even implies that $\cM(\Delta) \to
\cM(\Lambda)$ is surjective for any subfans $\Lambda \subset
\Delta\subset \Sig$.

Let $\Pure(\cA_\Sig)$ denote the full subcategory of $\cA\mof$
consisting of all pure sheaves, and let $K^b(\Pure(\cA))$ be the
category of bounded complexes of pure sheaves, with morphisms taken up
to chain homotopy.  Our main theorem is the following.

\begin{thm} \label{homotopy of pure sheaves} 
  The natural functor $K^b(\Pure(\cA)) \to D^b(\cA)$ is an equivalence
  of categories.
\end{thm} 

This sort of theorem is familiar when instead of pure objects we have
complexes of projective or injective objects in an abelian category.
The idea of the theorem is that a pure sheaf $\cM$ is half injective
and half projective: the locally free condition says that the stalk
$\cM(\sig)$ is a projective $\cA_\sig$-module, and flabbiness means
that $\cM$ is injective as a sheaf of vector spaces.

\subsection{Flabby sheaves}
Let us make more precise in what sense flabby sheaves are injective.
We first need the following definition.

\begin{defn} \label{strongly injective} 
  An injective map $\cM\to \cN$ is called {\em strongly injective} if
  the inclusion of stalks $\cM(\sig)\to \cN(\sig)$ splits for every
  $\sig\in \Sig$.
\end{defn}

\begin{prop} \label{quasi-injectivity}
  Suppose that $\cI$ is a flabby $\cA$-module.  If $\eta\colon\cM \to \cN$ is a
  strongly injective map, and $\cN$ is locally free, then the induced
  homomorphism $\Hom_{\cA}(\cN, \cI) \to \Hom_{\cA}(\cM, \cI)$ is
  surjective.
\end{prop}

\begin{proof} Take a map $\phi\colon\cM \to \cI$.
  We define a lift $\psi\colon \cN\to \cI$ of $\phi$ inductively on an
  increasing sequence of subfans.  Defining $\psi$ on the zero cone is
  trivial.  So suppose $\Delta\subset \Sig$ is a subfan with more than
  one cone, that $\tau\in \Delta$ is a maximal cone, and that
  $\psi|_{\Delta\setminus\{\tau\}}$ has been defined already.
  
  Since $\eta$ is strongly injective, we can choose a
  splitting of $\cA_\tau$-modules 
%\begin{equation} \label{splitting}
  $\cN(\tau) = \eta(\cM(\tau)) \oplus M$;
%\end{equation}
  since $\cN(\tau)$ is a free $\cA_\tau$-module, so are $\cM(\tau)$
  and $M$.  Define the restriction of $\psi$ to $\eta(\cM(\tau))$
  to be $\phi \eta^{-1}$.  To define $\psi$ on
  $M$, we need a map $M \to \cI(\tau)$ making the square
\[\xymatrix{
  M \ar[r]\ar[d] & \cI(\tau)\ar[d] \\
  \cN({\bdy\tau})\ar[r] & \cI({\bdy\tau}) }\] commute.  The right-hand
vertical map is surjective, since $\cI$ is flabby, and since $M$ is free, the
required map exists.
\end{proof}

\begin{prop} \label{fully faithful}
  If $\phi\colon \cM^\udot\to \cN^\udot$ is a quasi-isomorphism of
  pure complexes, then it has a homotopy inverse.
\end{prop}

\begin{lemma} 
\label{map from acyclic}
If $\cZ^\udot$ and $\cM^\udot$ are bounded complexes of pure sheaves,
and $\cZ^\udot$ is acyclic, then any map $\cZ^\udot \to \cM^\udot$ is
chain-homotopic to zero.
\end{lemma}

\begin{proof}
  Proposition \ref{quasi-injectivity} allows us to copy the standard
  argument used when the objects $\cZ^i$ are injective, provided that
  we know that each $\coker d_\cZ^i \to \cZ^{i+2}$ is strongly
  injective.  In other words, we need to show that $\coker
  (d_{\cZ^i(\sig)}) \to \cZ^{i+2}(\sig)$ is a split injection for all
  $\sig\in\Sig$.  This follows from the fact that $\cZ^\udot(\sig)$ is
  an acyclic complex of free $\cA_\sig$-modules.
\end{proof}

\begin{proof}[Proof of Proposition \ref{fully faithful}]
  Let $\cZ^\udot$ be the mapping cone of $\phi$.  Applying the lemma
  to the connecting map $\cZ^\udot\to \cM^\udot[1]$ gives a chain
  homotopy whose components are maps $h^i\colon \cZ^i \to \cM[1]^{i-1}
  = \cM^i$.  But $\cZ^i = \cN^i \oplus \cM^{i+1}$, so the first
  component of $h^i$ gives a map $\psi^i\colon \cN^i \to \cM^i$.  The
  resulting map of complexes $\psi\colon \cN^\udot \to \cM^\udot$ is a
  homotopy inverse of $\phi$.
\end{proof}

\subsection{Locally free sheaves}
Locally free sheaves act like projective objects in a very similar
way.  We do not need the following result, but we include it to
illustrate the parallels with the situation for flabby sheaves.

\begin{defn} A surjective morphism $\cM \to \cN$ between 
  $\cA_\Sig$-modules is called {\em strongly surjective} if the
  induced homomorphism $\cM(\sig, \bdy\sig) \to \cN(\sig,\bdy\sig)$ is
  surjective for every $\sig\in \Sig$.
\end{defn}

\begin{prop} \label{quasi-projectivity} Let $\cP$ be a locally free sheaf.  If
  $\cM \to \cN$ is strongly surjective and $\cM$ is flabby, then the
  homomorphism $\Hom_{\cA}(\cP,\cM) \to \Hom_{\cA}(\cP,\cN)$ is
  surjective.
\end{prop}

The proof is left to the reader.

\subsection{There are enough pure sheaves}  
To finish the proof of Theorem \ref{homotopy of pure sheaves} we need
to show that there are ``enough'' pure sheaves to represent any
complex in $D^b(\cA)$.  This follows from a two-step resolution
process, using the following result.

\begin{prop}  \label{enough pures} Take any object $\cM\in \cA\mof$.
\begin{enumerate}
\item[(a)] There exists a locally free $\cA$-module $\cP$ and a strong
  surjection $\cP \to \cM$.  If $\cM$ is flabby, then $\cP$ can be
  chosen to be pure.
\item[(b)] There exists a flabby $\cA$-module $\cI$ and a strong
  injection $\cM \to \cI$.  If $\cM$ is locally free, then $\cI$ can
  be chosen to be pure.
\item[(c)] If $\cM$ is zero on the subfan $\bdy\sig$ for some
  $\sig\in\Sig$, then the maps in {\rm (a)} and {\rm (b)} can be
  chosen to be isomorphisms on all of $[\sig]$.
\end{enumerate} 
\end{prop}

Assuming this result, we can now prove Theorem \ref{homotopy of pure
  sheaves}.  Consider any complex $\cM^\udot \in D^b(\cA)$.  Part (a)
of the proposition allows us to find a complex $\cP^\udot$ of locally
free sheaves and a quasi-isomorphism $\cP^\udot \stackrel{\sim}{\to}
\cM^\udot$.  Note that (c) implies that $\cP^\udot$ can be chosen to
be a bounded complex.

Part (b) then implies that there is a complex $\cI^\udot$ of pure
sheaves and a quasi-isomorphism $\cP^\udot \stackrel\sim\to
\cI^\udot$.  Note that here it is crucial that (b) gives strong
injections: to construct $\cI^j$ we apply (b) to $\cM =
\coker(\cP^{j-1} \to \cP^j \oplus \cI^{j-1})$, which is locally free
since $\cP^j$ and $\cI^{j-1}$ are and $\cP^{j-1} \to \cI^{j-1}$ is a
strong injection.  Using (c) again, we see that $\cI^\udot$ can be
chosen to be a bounded complex.

Thus the functor $K^b(\Pure(\cA)) \to D^b(\cA)$ is essentially
surjective.  The same two-step resolution process also shows it is
fully faithful, by the following well-known result.  Let $\cC$ be an
abelian category, $K(\cC)$ the homotopy category of complexes in
$\cC$.
\begin{lemma}  Let $K_1 \subset K_2$ be full triangulated subcategories
  of $K(\cC)$, and let $D_1$, $D_2$ be the corresponding derived
  categories.  If either of the following conditions holds, then the
  functor $D_1 \to D_2$ is fully faithful.
\begin{enumerate}
\item For any quasi-isomorphism $X^\udot \to Y^\udot$ in $K_2$, with
  $X^\udot$ in $K_1$, there exists a quasi-isomorphism $A^\udot \to
  X^\udot$ with $A^\udot \in K_1$.
\item For any quasi-isomorphism $X^\udot \to Y^\udot$ in $K_2$, with
  $Y^\udot$ in $K_1$, there exists a quasi-isomorphism $Y^\udot \to
  B^\udot$ with $B^\udot \in K_1$.
\end{enumerate}
\end{lemma}

\begin{cor} \label{loc free --> pure}
  Suppose $\cM$, $\cP$ are $\cA$-modules, $\cM$ is locally free, and
  $\cP$ is pure.  Then $\Hom_{D^b(\cA)}(\cM,\cP[i]) = 0$ unless $i =
  0$.
\end{cor}

\begin{proof} We can replace $\cM$ by a 
  quasi-isomorphic complex $\cI^\udot$ of pure sheaves, with $\cI^j =
  0$ for $j < 0$ and $\coker(\partial^{j-1}) \to \cI^{j+1}$ strongly
  injective for all $j \ge 0$.  Now apply Proposition
  \ref{quasi-injectivity}.
\end{proof}

\begin{proof}[Proof of Proposition \ref{enough pures}]
  To prove (a), we construct the object $\cP$ and the map $\phi\colon
  \cP \to \cM$ simultaneously, by induction on subfans.  For the base
  case when $\Delta = \{o\}$, we set $\cI(o) = \cM(o)$ and let
  $\phi|_o$ be the identity map.
  
  Now suppose $\Delta$ is a fan with at least two cones, $\tau$ is a
  maximal cone, and the restrictions of $\cP$ and $\phi$ to $\Delta
  \setminus \{\tau\}$ have already been defined.  To define them on
  all of $\Delta$ it is enough to define them on $[\tau]$, since the
  resulting sheaves and maps can be glued.
  
  This amounts to choosing a free $\cA_\tau$-module $\cP(\tau)$ and
  homomorphisms $\bdy_\cP\colon \cP(\tau) \to \cP(\bdy\tau)$ and
  $\phi_\tau \colon \cP(\tau) \to \cM(\tau)$ so that $\phi_\tau$ and
  the induced map $\ker \bdy_\cP \to \ker \bdy_\cM$ are surjective and
  the square
\[\xymatrix{
  \cP(\tau) \ar[r]^{\phi_\tau} \ar[d]_{\bdy_\cP} & \cM(\tau)\ar[d]^{\bdy_\cM}\\
  \cP({\bdy\tau}) \ar[r]_{\phi|_{\bdy\tau}} & \cM({\bdy\tau}) }\]
commutes. To do this, find free $\cA_\tau$-modules $M_1$,$M_2$ and
maps $p_1\colon M_1 \to \cP({\bdy\tau})$ and $p_2\colon M_2 \to
\cM(\tau)$ so that $\im p_1 = (\phi|_{\bdy\tau})^{-1}(\im \bdy_\cM)$
and $\im p_2 = \ker(\bdy_\cM)$.  Then let $\cP(\sig) = M_1 \oplus
M_2$, and $\bdy_\cP = p_1\oplus 0$.  To define $\phi_\tau$, let
$\phi_\tau|_{M_2} = p_2$, and for each $a$ in a basis of $M_1$, define
$\phi_\tau(a)$ to satisfy $(\phi|_{\bdy\tau})(p_1(a)) =
\bdy_\cM\phi_\tau(a)$.

To prove (b), we again proceed by induction.  The base case is again
trivial, and we are reduced to the problem of extending the sheaf
$\cI$ and morphism $\phi\colon \cM\to \cI$ from $\bdy\tau$ to $[\tau]$
as before.  This, in turn, amounts to finding an $\cA_\tau$-module
$\cI(\tau)$, a surjective restriction homomorphism $\bdy_\cI\colon
\cI(\tau) \to \cI({\bdy\tau})$, and a split injection $\phi_\tau\colon
\cM(\tau) \to \cI(\tau)$, such that the square
\[\xymatrix{
  \cM(\tau) \ar[r]^{\phi_\tau} \ar[d]_{\bdy_\cM} & \cI(\tau)\ar[d]^{\bdy_\cI}\\
  \cM({\bdy\tau}) \ar[r]_{\phi|_{\bdy\tau}} & \cI({\bdy\tau}) }\]
commutes.  This can be done by letting $\cI(\tau) = \cI({\bdy\tau})
\oplus \cM(\tau)$, and letting $\phi_\tau = (0, \id_\cM)$ and
$\bdy_\cI = \id_{\cI({\bdy\tau})} \oplus (\phi|_{\bdy\tau} \circ
\bdy_\cM)$.

For the second statement of (b), we make a different choice at the
inductive step: take a free $\cA_\tau$-module $M$ and a
surjective homomorphism $p\colon M\to \cI({\bdy\tau})$.  We then
define $\cI(\tau) = M \oplus \cM(\tau)$ and $\bdy_\cI = p \oplus 0$.
Since $\cM(\tau)$ is free by assumption, $\cI(\tau)$ is free as well.
The required map $\phi_\tau$ now exists because $\bdy_\cI$ is
surjective and $\cM(\tau)$ is free.

Checking that these constructions satisfy (c) is easy.
\end{proof}

\subsection{Indecomposible pure sheaves}
The notion of pure $\cA$-module was first used in \cite{BrLu,BBFK} to
model direct sums of (shifted) intersection cohomology sheaves.  The
indecomposible pure objects are models of single intersection
cohomology sheaves.  We recall here their basic properties.

For a cone $\sig\in\Sig$, let $c(\sig)$ denote the codimension of
$\sig$ in the ambient vector space.
For $n\in \Z$, recall that $\lb n\rb \colon D^b(\cA) \to D^b(\cA)$ is the 
functor which shifts the degree {\em down} by $n$.

\begin{thm} \label{pure simples}
  For every $\sig\in\Sig$, there is an indecomposible pure
  $\cA$-module $\cL^\sig$, unique up to a scalar isomorphism, for
  which (1) $\cL^\sig(\tau) = 0$ unless $\sig \prec\tau$, and (2)
  $\cL^\sig(\sig) = \cA_\sig\lb c(\sig)\rb $.
  
  These objects satisfy the following:
\begin{enumerate} 
\item Every pure sheaf is isomorphic to a finite direct sum $\oplus_i
  \cL^{\sig_i}\lb n_i\rb $ with $\sig_i\in \Sig$ and $n_i \in \Z$.
\item For all $\tau \in \Sig\setminus\{\sig\}$, $\cL^\sig(\tau)$ is
  generated in degrees $< -c(\tau)$.
\item For all $\tau \in \Sig$, $\cL^\sig({\tau,\bdy\tau})$ is a free
  $\cA_\tau$-module; it is generated in degrees $> -c(\tau)$, unless
  $\sigma = \tau$.
\end{enumerate}
\end{thm}
\begin{rmk}
  We put the generator of $\cL^\sig(\sig)$ in degree $-c(\sig)$
  (rather than degree $0$ as in \cite{BrLu,BBFK}) so that the
  resulting object will be perverse, i.e.\ in the core of the
  $t$-structure which we define in the next section.
\end{rmk}

A proof of (1) can be found in \cite{BrLu,BBFK}, while (2) and (3)
follow from work of Karu \cite{Ka}.

Next we look more carefully at homomorphisms between pure sheaves.
Because of Theorem \ref{pure simples}, it is enough to look at the
objects $\cL^\sig\lb n\rb $.
\begin{thm} \label{Homs between simples}
  Take $\sig,\tau\in \Sig$, $n\in \Z$.  Let $H =
  \Hom_{\cA}(\cL^\sig,\cL^\tau\lb n\rb )$.
\begin{enumerate}
\item[(a)] If $n < 0$, then $H=0$.
\item[(b)] If $n=0$, then $H=0$ unless $\sig = \tau$, in which case
  $\dim_\R H = 1$, with a basis given by the identity map.
\item[(c)] If $n = 1$ and $\sig\prec\tau$, then restricting to $\tau$
  gives an isomorphism
\[ H \cong 
\Hom_{\cA_\tau}(\cL^\sig(\tau),\cL^\tau({\tau,\bdy\tau}) ) =
\cL^\sig(\tau)^*_{c(\tau)-1},\] while if $\tau\prec\sig$, restricting
to $\sig$ gives an isomorphism
\[H \cong 
\Hom_{\cA_\sig}(\cL^\sig(\sig),\cL^\tau(\sig,\bdy\sig)) =
\cL^\tau(\sig,\partial\sig)_{-c(\sig)+1}.\] If $\tau$ and $\sig$ are
not comparable, then $\Hom_{\cA}(\cL^\sig,\cL^\tau\lb 1\rb ) = 0$.
\end{enumerate}
\end{thm}

The proof is by induction, using Theorem \ref{pure simples} and the
following lemma.

\begin{lemma} Suppose $\cM,\cN $  are $\cA$-modules, $\cM$ 
  is locally free, and $\cN$ is flabby.  If $\sig\in \Sig$ is a
  maximal cone, then there is a short exact sequence
\[0 \to \Hom_{\cA_\sig}(\cM(\sig),\cN(\sig,\bdy\sig))
\to \Hom_{\cA}(\cM,\cN)\]\[ \to
\Hom_{\cA|_{\Sig\setminus\{\sig\}}}(\cM|_{\Sig\setminus\{\sig\}},
\cN|_{\Sig\setminus\{\sig\}}) \to 0.\]
\end{lemma}

The following corollary of Theorem \ref{Homs between simples}(b) is
useful.
\begin{cor} \label{Cor to Homs between simples}
  If $\cM \to \cN$ is a morphism between two pure sheaves, each of
  which is a direct sum of various $\cL^\sig$ (without degree shifts),
  then the kernel and cokernel are both pure.
\end{cor}

We will also need the following technical lemma.
\begin{lemma} \label{Noetherian}
  Suppose that $\cM$ is a pure $\cA$-module.  Then the graded
  endomorphism ring
\[R = \End_{\cA\Mod}(\cM) = \oplus_{n\in \Z}\Hom_{\cA\mod}(\cM,\cM\lb n\rb )\]
is Noetherian.
\end{lemma}
\begin{proof}
  There is a homomorphism from $\cA_o$ (polynomial functions on
  $N\otimes \C$) to $R$ given by pointwise multiplication.  The ring
  $R$ is contained in $\oplus_{\tau \in \Sig}
  \End_{\cA_o\Mod}(\cM(\tau))$, which is a finitely generated
  $\cA_o$-module.
\end{proof}

\subsection{Perverse $t$-structure} \label{t-structure on cA-complexes}
We define a $t$-structure on the triangulated category $D^b(\cA)$,
analogous to the usual one on the equivariant derived category
$D^b_T(X)$ whose core consists of equivariant perverse sheaves.

Let $K^{\ge 0}$ (respectively $K^{\le 0}$) be the full subcategory of
$K^b(\Pure(\cA))$ consisting of complexes which are quasi-isomorphic
to a complex $\cM^\udot$, where $\cM^i \cong \bigoplus_k
\cL^{\sig_k}\lb n_k\rb $, $\sig_k\in \Sig$, $n_k \le i$ (respectively $n_k
\ge i$).

\begin{thm} \label{t-structure thm}
  This defines a $t$-structure on $K^b(\Pure(\cA))$.  The heart
  $P(\Sig) = K^{\ge 0}\cap K^{\le 0}$ is equivalent to the full
  subcategory of $P(\Sig)$ consisting of bounded complexes $\cP^\udot$
  so that for any $i$ we have
\begin{equation*} \tag{*} \cP^i \cong
\bigoplus_{k=1}^l \cL^{\sig_k}\lb i\rb 
\end{equation*}
with $\sig_1,\dots,\sig_l \in \Sig$.
\end{thm}   

Since $D^b(\cA)$ and $K^b(\Pure(\cA))$ are equivalent categories, this
defines a $t$-structure $(D^{\le 0}, D^{\ge 0})$ on $D^b(\cA)$ as
well.

\begin{rmk} Note that all chain homotopies between complexes satisfying
  (*) automatically vanish, by Theorem \ref{Homs between simples}.
  Thus if objects in $P(\Sig)$ are represented by such complexes,
  morphisms are just morphisms of complexes.

  The resulting category of mixed equivariant perverse sheaves is similar
to a construction of Vybornov \cite{V}. 
\end{rmk}

\begin{proof}
  There are four things to check to show that $(K^{\ge 0}, K^{\le 0})$
  forms a $t$-structure.  It is clear that $K^{\ge 0} \subset K^{\ge
    0}[1]$ and $K^{\le 0}[1] \subset K^{\le 0}$.  If $\cM^\udot \in
  K^{\le 0}$ and $\cN^\udot \in K^{\ge 1} = K^{\ge 0}[-1]$, we have
  $\Hom(\cM^\udot,\cN^\udot) = 0$, by Theorem \ref{Homs between
    simples}(a).  Given a distinguished triangle
\begin{equation*}
\cE^\udot \to \cM^\udot \to \cN^\udot \stackrel{[1]}{\to},
\end{equation*}
where $\cE^\udot$ and $\cN^\udot$ are both in $K^{\le 0}$ (resp.\ 
$K^{\ge 0}$), then $\cM^\udot \in K^{\le 0}$ (resp. $\cM^\udot\in
K^{\ge 0}$), since the triangle comes from a short exact sequence $0
\to \cE^\udot \to\wt\cM^\udot \to \cN^\udot \to 0$, with $\wt\cM^\udot
\cong \cM^\udot$.

Finally, we need to show that for any $\cM^\udot\in K^b(\Pure(\cA_\Sig))$
there exists a triangle $\cE^\udot \to \cM^\udot \to \cN^\udot
\stackrel{[1]}{\to}$ with $\cE^\udot \in K^{\le 0}$ and $\cN^\udot \in
K^{\ge 1}$.  To do this, we write each $\cM^i$ as a direct sum
$\oplus_j \cM^i_j$, where $\cM^i_j$ is isomorphic to a sum of various
$\cL^\sig\lb j\rb $.  Then we can write the differential $d^i\colon \cM^i
\to \cM^{i+1}$ as a sum $\sum\limits_{\stackrel{\scriptstyle j,k\in
    \Z}{k \ge 0}}d^i_{jk}$, where $d^i_{jk}\colon \cM^i_j \to
\cM^{i+1}_{j+k}$.

We then let
\[\cE^i =  \ker d^i_{i,0} \oplus\bigoplus_{i-j < 0} \cM^i_j ,\]
which is pure by Corollary \ref{Cor to Homs between simples}.  It is a
subcomplex of $\cM^\udot$, and it clearly lies in $K^{\le 0}$.
Moreover, it is compatible with the decomposition $\cM^i = \oplus_j
\cM^i_j$, so if we let $\cN^\udot = \cM^\udot/\cE^\udot$, we have a
decomposition $\cN^i = \oplus \cN^i_j$ compatible with the quotient
map.  Let $\tilde d$ be the differential of $\cN^\udot$; it can be
decomposed $\tilde d^i = \oplus \tilde d^i_{jk}$ as before.

We must show that $\cN^\udot$ lies in $K^{\ge 1}$.  Note that letting
\[\cN^i = \wt\cN^i_{i} \oplus \im(\tilde{d}^i_{i,0})\]
defines a subcomplex $\cN^\udot$ of $\cN^\udot$ which is
quasi-isomorphic to $0$.  Since $\cN^\udot/\cN^\udot$ is clearly in
$K^{\ge 1}$, so is $\cN^\udot$.

For the second statement, suppose $\cM^\udot$ is in $K^{\le 0}\cap
K^{\ge 0}$.  Then $\wt{\cM}^i = (\ker d^i_{i,0})/(\im
d^{i-1}_{i-1,0})$ gives a quasi-isomorphic complex which satisfies
(*).
\end{proof}

\subsection{$t$-exactness}
In this section, we prove
\begin{thm} \label{t-exact}
The functor $F_T\colon D^b(\cA) \to D^b_T(X_\Sig)$
  defined in \S\ref{equivariant sheaves and DG modules} is $t$-exact.
\end{thm} 
Here $D^b(\cA)$ has the
$t$-structure just defined, and $D^b_T(X)$ has the perverse
$t$-structure from \cite{BL}.

In terms of the presentation we use of $D^b_T(X)$ as a full
subcategory of $D^b(X_T)$, we define this second $t$-structure in
terms of the usual perverse $t$-structure $(D^{\le 0}(X), D^{\ge
  0}(X))$ on $D^b(X)$.  Since $X_T$ is a fiber bundle over $BT$ with
fiber $X$, we get an embedding of $X$ into $X_T$ by choosing a
basepoint in $BT$.  This gives rise to a ``forgetful functor''
\[\For\colon D^b_T(X) \to D^b(X)\]
which is simply restriction to $X$.  Then our perverse $t$-structure
is $(\For^{-1}D^{\le 0}(X),\For^{-1}D^{\ge 0}(X))$.

For a face $\sig\in \Sig$ let $j_\sig$ be the
inclusion of $O_\sig \hookrightarrow X$.  Let $F_{T,\Sig}\colon
D^b(\cA) \to D^b_T(X)$ and $F_{T,\sig}\colon D^b(\cA_\sig\mof)\to
D^b_T(O_\sig)$ be the realization functors.  We have restriction
functors
\[j_\sig^*,j^!_\sig\colon D^b_T(X) \to D^b_T(O_\sig)\]
which are simply restriction and corestriction to $O_{\sig,T}$ (note
that although the space $X_T$ is not locally compact, the
corestriction can be defined using the derived ``restriction with
supports'' functor $\bR\Gamma_{O_{\sig,T}}$; all the usual adjunction
and base change properties still apply).

Theorem \ref{t-exact} follows from Theorem \ref{pure simples} and the
following result, which describes the stalk and costalk functors
in terms of $\cA$-modules.  The proof will be given in \S\ref{appendix}.

\begin{thm}  If $\cM^\udot \in K^b(\Pure(\cA))$, there are natural 
isomorphisms in $D^b_T(O_\sig)$:
\begin{itemize}
\item[(a)] $j_\sig^*F_{T,\Sig} \cM^\udot \simeq
  F_{T,\sig}(\cM^\udot(\sig))$,
\item[(b)] $j_\sig^!F_{T,\Sig} \cM^\udot \simeq
  F_{T,\sig}(\cM^\udot(\sig,\partial \sig))$.
\end{itemize}
Let $i_x$ be the inclusion of a point $x$ into $ O_{\sig,T}$.  Then 
there are natural isomorphisms in $D^b(\C\mof) = D^b(pt)$
\begin{itemize}
\item[(c)] $i_x^*j_\sig^*F_{T,\Sig} \cM^\udot \simeq \C
  \otimes_{\cA_\sig} \nu\cM^\udot(\sig)$,
\item[(d)] $i_x^*j_\sig^!F_{T,\Sig} \cM^\udot \simeq \C
  \otimes_{\cA_\sig} \nu\cM^\udot(\sig,\partial \sig)$.
\end{itemize}
\end{thm} 

Here $\nu$ is the functor $D(\cA_\sig\mod) \to
D(\text{DG-$\cA_\sig$})$ of Example \ref{forget the grading}.
We can make a similar statement for general complexes $\cM^\udot 
\in D^b(\cA)$ if we replace the functors on the right hand
sides of (a), (b), (c), and (d) with the appropriate derived
functors: replace $\otimes_{\cA_\sig}$ by $\stackrel{L}{\otimes}_{\cA_\sig}$
and $-(\sig,\partial \sig) = \Gamma_\sig(-)|_{\{\sig\}}$ by 
$\bR\Gamma_\sig(-)|_{\{\sig\}}$.

\subsection{Mixed structure} \label{mixed equivariant sheaves}

For any $n\in \Z$, 
define an automorphism of $D^b(\cA)$ by $\la n\ra = [n]\lb -n\rb $.  
It is obviously $t$-exact, so it induces an automorphism of
$P(\Sig)$.  Define a functorial filtration on objects of
$P(\Sig)$ by $W_j\cP^\udot = \oplus_{i\ge -j} \cP^i$, assuming
$\cP^\udot$ is a complex satisfying condition (*) of Theorem
\ref{t-structure thm}.  It is easy to see this defines a mixed structure on
$P(\Sig)$.  We will show in the next section that the functor
$F_T\colon P(\Sig)\to P_T(X_\Sig)$ is a grading on $P_T(X_\Sig)$
in the sense of \S\ref{grading section}.

\section{The toric Koszul functor} \label{main proofs}
\subsection{} 
Let $X = X_\sig$ be a normal affine $T$-toric variety with a single
fixed point defined by a
cone $\sigma \subset N_T\otimes \bbR$,
with $\dim \sig = \rank N_T = n$. Let $\check X$ be
the dual $\check T$-toric variety defined by the dual cone $\check
\sigma \subset N_{\check T}\otimes \bbR$. 
Put $\Sig = \Sig _X =[\sigma]$, 
$\check \Sig = \Sig _{\check X}=[\check \sigma]$,  $\cA = \cA_\Sig$,
and $\cT = \cT_{\check \Sig}$.  Recall the identification 
of ringed quivers
\begin{equation}\label{aaa}
(\Sig,\cA) = ((\check\Sig)^\circ,\cT)
\end{equation}
from \S\ref{dual toric varieties}.

In this section we define our Koszul equivalence
\[K\colon D^b(\cA\mof)\to D^b(LC_\cF(X^\vee)).\]
It will be a composition of three equivalences:
\[D^b(\cA\mof) \stackrel{\kappa}{\to} D^b(\co\cA\mcf)
\to D^b(LC_\cF(X^\vee)) \stackrel{\la -n\ra}\to D^b(LC_\cF(X^\vee));\]
The last functor $\la -n \ra$ is the twist defined in \S\ref{F-sheaves}.
The middle functor is the equivalence
of Theorem \ref{combinatorial F-sheaves} combined with 
\eqref{aaa}.  The functor $\kappa$ is a combinatorial form of Koszul duality
which makes sense for any fan, rational or not.  We define it in the 
next section.

Our main results can be summarized as follows.
\begin{thm} \label{K is Koszul}
$K$ is a Koszul equivalence in the sense of Definition
\ref{Koszul functor definition}.  Here we use the $t$-structures 
and mixed structure on 
$D^b(\cA\mof)$ and $D^b(LC_\cF(X^\vee))$ defined in \S\ref{mixed F-sheaves} and
\S\ref{mixed equivariant sheaves}, and the ring $R$, resp.\ $R^\vee$, is
the opposed ring of the 
graded endomorphism ring of a mixed projective generator of $P(\cA_\Sig)$ (resp.\  
a mixed injective generator of $P_\cF(X^\vee)$). 
\end{thm}
  
\subsection{Combinatorial Koszul functor} \label{combinatorial Koszul}
Fix a fan $\Sig$ (rational or not) in $\bbR ^n$ with
the corresponding ``structure sheaf" $\cA =\cA _\Sig$.
We define the functor $\kappa\colon D^b(\cA\mof) {\to} D^b(\co\cA\mcf)$
as follows.

%There is a ``naive" anti-equivalence $(\cdot )^*$ between
%these categories. Namely, let $A$ be a graded polynomial
%ring and $M$ be a finitely generated graded $A$-module. Consider
%its graded dual $M^*$, where $(M^*)_n=\Hom (M_{-n},\bbC)$. Then
%$M^*$ is a co-finite $A$-module and this gives the
%anti-equivalence
%$$(\cdot )^*\colon \cA _f\text{-Mod}\to \co\cA _{cf}\text{-Mod}.$$
%
%It follows that $(\cdot )^*$ maps projective objects in $\cA _f\text{-Mod}$
%to injectives in $\co\cA _{cf}\text{-Mod}.$ There are
%enough projectives in  $\cA _f\text{-Mod}$: for each $\tau \in \Sig$ and
%$\cM \in  \cA _f\text{-Mod}$ denote by $\cM _{[\tau]}$ the extension by
%zero to $\Sig$ of the restriction $\cM \vert _{[\tau]}$; then $\cA
%_{[\tau]}$ is projective. Hence $\cA^*_{[\tau]}$ is injective.
%
%This naive equivalence is contravariant, while we want a covariant
%equivalence.  We do this as follows.

\begin{defn}
  Fix an abelian category $\cC$.
\begin{enumerate}
  
\item[a)] A $\Sig$-{\it diagram} in $\cC$ is a collection $\{ \cM
  _\tau\} _{\tau\in\Sig}$ of objects of $\cC$ together with
  with morphisms $p_{\tau \xi}\colon \cM _\tau \to \cM _\xi$ for $\tau
  \geq \xi$, satisfying $p_{\rho\xi}p_{\tau\rho} = p_{\tau\xi}$
  whenever $\tau \geq \rho \geq \xi$.
  
\item[b)] Fix an orientation of each cone in $\Sig$. Then every $\Sig
  $-diagram $\cM=\{ \cM _\tau \}$ gives rise to the corresponding {\it
    cellular complex} in $\cC$:
  $$C^\bullet(\cM)=\ \bigoplus_{\dim (\tau)=n}\cM _\tau \to
  \bigoplus_{\dim (\xi)=n-1}\cM _\xi \to ...$$
  where the terms
  $\cM_\rho$ appear in degree $- \dim \rho$, and the differential is
  the sum of the maps $p_{\tau \xi}$ with $\pm$ sign depending on
  whether the orientations of $\tau $ and $\xi$ agree or not.
\end{enumerate}
\end{defn}

\begin{lemma} \label{Cech complex of constant sheaf}
  Let $\cM=\{\cM _\tau \}$ be a constant $\Sig$-diagram supported
  between cones $\eta $ and $\xi$ in $\Sig$. That is $\cM _\tau =M$
  for a fixed $M$ if $\eta \geq \tau \geq \xi$, and $\cM_\tau =0$
  otherwise; for $\eta \geq \tau _1 \geq \tau _2 \geq \xi$ the maps
  $p_{\tau _1\tau _2}$ are the identity.  If $\eta\ne \xi$, then the
  cellular complex $C^\bullet (\cM)$ is acyclic.
\end{lemma}

\begin{proof} The complex $C^\bullet (\cM)$ is isomorphic to an 
  augmented cellular
  chain complex of a closed ball of dimension $\dim (\eta )-\dim
  (\xi)-1$.
\end{proof}

Recall that the sheaves $\cA ^*_{[\tau ]}$ are injective objects of
co-$\cA \text{-Mod}$ for every $\tau \in \Sig$.  Consider the
$\Sig$-diagram $\cK=\{ \cK _\tau \}$ in co-$\cA \text{-Mod}$, where
$\cK _\tau =\cA ^*_{[\tau ]}$ and the maps $p_{\tau \xi}$ are the
projections.  This diagram $\cK$ defines a covariant functor
$$\kappa \colon D^b(\cA\mof)\to D^b(\co\cA\mcf)$$
in the following
way. If $\cN \in \cA\mof$ is locally free, the collection \[\cK
\otimes \cN=\{ \cK _\tau \otimes _{\cA _\tau}\cN _\tau \}\] is a $\Sig
$-diagram of co-$\cA$-modules.  Thus its cellular complex $C^\bullet
(\cK \otimes \cN)$ is a complex of co-$\cA$-modules.  By Theorem
\ref{homotopy of pure sheaves}, every
object in $D^b(\cA\mof)$ is quasi-isomorphic to a complex of locally free
(in fact, pure) $\cA$-modules, so we obtain a derived functor
$$\kappa(\cdot )=C^\bullet(\cK \stackrel{\bbL}{\otimes}\cdot)\colon
D^b(\cA\mof)\to D^b(\co\cA\mcf).$$
We call it the combinatorial Koszul
functor.  

\begin{lemma} \label{calculations of kappa}
  For any cone $\tau\in \Sig$, there are isomorphisms $\kappa
  (\cA_{\{\tau \}})\cong \cA^*_{[\tau]}[\dim\tau]$ and $\kappa
  (\cA_{[\tau ]})\cong\cA^*_{\{ \tau\}}[\dim\tau]$.
\end{lemma}

The first isomorphism is obvious; the second follows from Lemma
\ref{Cech complex of constant sheaf}.

\begin{prop} The functor $\kappa$ is an equivalence of triangulated categories.
\end{prop}

\begin{proof} 
  The category $D^b(\cA _f\text{-Mod})$ is the triangulated envelope of
  either all objects of the form $\cA_{[\tau]}\lb k\rb$ 
  or all objects of the form $\cA _{\{\tau\}} \lb k\rb$, in either case taken over all 
  $\tau \in \Sig$ and $k\in \Z$.  Similarly, $D^b(\co\cA_{cf}\text{-Mod})$) is the 
  triangulated envelope of all objects of the form $\cA^*_{[\tau]} \lb k\rb$ or
  all objects of the form $\cA^*_{\{\tau \}} \lb k\rb$.
  
  So it suffices to show that for any $k,l,\tau , \xi$ the functor
  $\kappa$ induces an isomorphism
  $$\kappa \colon \Hom _{\cA _f\text{-Mod}}(\cA _{[\tau]}\lb k\rb , \cA
  _{\{\xi \}}\lb l\rb )\to \Hom _{\co\cA _{cf}\text{-Mod}}(\cA ^*_{\{\tau
    \}}\lb k\rb ,\cA ^*_{[\xi]}\lb l\rb ).$$
  Both sides are equal to the $l-k$
  graded part of $\cA _{\tau}$ if $\tau =\xi$ and vanish otherwise.
\end{proof}
Thus $K$ satisfies property (1) of the definition of
a Koszul equivalence (Definition 
\ref{Koszul functor definition}).   Property (2) follows immediately
from the definition of the twist functors $\la n\ra$ in the 
categories $D^b(\cA)$ and $D^b(LC_\cF(\check X))$.
Showing that $K$ sends simples to injectives and indecomposable
projectives to simples (properties (3) and (4)) will take up the
remainder of the paper.

\begin{rmk}  We can think of $\kappa$ as convolution with
kernel $\cK$.  This is more enlightening if one considers $\cK$ 
as a sheaf on $\Sig \times (\check \Sig)^\circ$ by the natural 
identification $\Sig = (\check \Sig)^\circ$.  
The support of $\cK$ is then the  ``combinatorial conormal variety'' 
\[\Lambda = \{(\tau,\alpha)\in \Sig \times (\check \Sig)^\circ \mid
 \tau^\bot \le \alpha\}.\]  If $p_1\colon \Lambda \to \Sig$, $p_2\colon \Lambda \to
(\check \Sig)^\circ$ are the projections,
then $\cK =  p_2^{-1}\cT^*$, using the identification \eqref{aaa}.  
Note that $\cK$ has a natural action of $p_1^{-1}\cA$ which 
commutes with the action of $p_2^{-1}\cT$.  In fact, 
$\Lambda$ is the largest subset of $\Sig \times (\check \Sig)^\circ$
for which this is true.  
\end{rmk}

\subsection{Proof of Theorem \ref{K is Koszul}, part I: 
  $K(\text{simple})$ is injective}

Let us examine what the functor $K$ does to indecomposable pure
sheaves.  For a face $\alpha \in [\check\sig]$, define $I^\udot_\alpha
= K(\cL^{\alpha^\bot})$.  It is perverse, as follows from the
following more general statement.  Recall the objects
$\nabla_\alpha^\udot \in P_\cF(\check X)$ from \S\ref{perverse
  t-structure}.
\begin{prop} \label{locally free to perverse}
  If $\cM^\udot\in D^b(\cA\mof)$ is given by placing a locally free
  $\cA$-module in degree $0$, then $K(\cM^\udot)$ is perverse, and has
  a filtration whose graded pieces are objects $\nabla_\alpha^\udot\la
  k\ra$, $\alpha \in \check \Sig$, $k \in \Z$.
\end{prop}

\begin{proof}  A locally free $\cA$-module has a filtration whose
  subquotients are sheaves $\cA_{\{\tau\}}\{ k\}$, $\tau\in \Sig$,
  $k\in \Z$. By Lemma \ref{calculations of kappa}, we have
  $K(\cA_{\{\tau\}}) \cong \nabla_{\tau^\bot}^\udot\la \dim \tau - n
  \ra$.  The result follows.
\end{proof}

\begin{thm} \label{K(simple) is injective}
  $I^\udot_\alpha$ is an injective object in $P_\cF(X)$, and its
image $F_{cf}(I^\udot_\alpha)$ is injective in $P_{cf}(X)$.

  With respect
  to the mixed structure defined in \S\ref{mixed F-sheaves}, 
we have $W_0I^\udot_\alpha \cong L^\udot_\alpha$.
\end{thm}

\begin{proof}
  To show that $I^\udot_\alpha$ is injective, we will show that the
  following statement holds for any $S^\udot \in P_\cF(\check X)$:
\begin{equation}
\tag{*} \Hom_{D^b(LC_\cF(\check X))}(S^\udot, I^\udot_\alpha[k]) =
  0 \;\text{for all}\; k > 0.
\end{equation}
First note that if
\[0\to S_1^\udot \to S_2^\udot \to S_3^\udot \to 0\]
is a short exact sequence in $P_\cF(\check X)$ and (*) holds for
$S_1^\udot$ and $S_3^\udot$ or for $S_2^\udot$ and $S_3^\udot$, then
it also holds for all three objects.

Note that (*) holds for $S^\udot = \nabla_\beta^\udot\la k\ra$ for any
$\beta\in [\check\sig]$ and $k\in \Z$, by applying $K$ to Corollary 
\ref{loc free --> pure}.  Thus (*) holds for any object of the form
$S^\udot = \bR j_{\beta *}E_\beta[\dim O_\beta]$, where $E_\beta$ is
any object in $LC_\cF(O_\beta)$, since $\bR j_{\beta *}$ is $t$-exact
(Proposition \ref{lower star is exact}), 
and $E_\beta$ can be resolved by a finite complex of injective
$\cF$-local systems, i.e.\ by direct sums of copies of objects
$\Theta_\beta\la k\ra$, $k\in \Z$.

Now we prove (*) for general $S^\udot$, by induction on the number of
orbits in the support.  If $O_\beta$ is an open orbit contained in
$\supp S^\udot$, then $S^\udot|_{O_\beta}$ is a $\cF$-local system
placed in degree $-\dim O_\beta$.  Consider the adjunction map
$\phi\colon S^\udot \to \bR j_{\beta *}(S^\udot|_{O_{\beta}})$, and
note that (*) holds for the target of $\phi$ by the previous
paragraph.  If $\supp S^\udot$ consists of the unique closed orbit
$O_{\sig^\vee}$, so $\beta = \check \sig$, then $\phi$ is an
isomorphism, and we are done.  Otherwise note that (*) holds by
induction for the kernel and cokernel of $\phi$, since they have
strictly smaller support, and thus it holds for $S^\udot$.  

Thus
$I^\udot_\alpha$ is injective.  The same argument shows that
$F_{cf}(I^\udot_\alpha)$ is injective in $P_{cf}(X)$. 

Apply Proposition \ref{locally free to perverse} to obtain a
filtration
\[M^\udot_0 \subset \dots \subset M^\udot_l = I^\udot_\alpha\]
with $M_0^\udot = \nabla^\udot_\alpha$ , and where the
$M^\udot_i/M^\udot_{i-1}$ for $i >0$ are sums of objects
$\nabla^\udot_\beta\la k\ra$ with $\beta \in [\alpha]\setminus
\{\alpha\}$ and $k \in \Z$.  In fact, using Lemma \ref{calculations of
  kappa} and property (2) of Theorem \ref{pure simples}, we see that
only twists $k > 0$ can occur.  The remaining statements of the
theorem follow using Corollary \ref{weights of nabla}.
\end{proof}

Let $I^\udot = \oplus_{\alpha \in [\sig^\vee]} I^\udot_\alpha$, so
$I^\udot = K(\cL)$, where $\cL = \oplus_{\tau\in [\sig]} \cL^\tau$.
\begin{prop} $I^\udot$ is a mixed injective generator (\S\ref{mixed categories}) 
of $P_\cF(\check X)$; $F_{cf}(I^\udot)$ is an
  injective generator of $P_{cf}(\check X)$.
\end{prop}
\begin{proof} 
  Let $Inj$ be the category of all finite direct sums of objects
of the form $I^\udot_\alpha\la k\ra$, $\alpha \in [\check \sig]$, $k\in
\Z$.  We need to show that any object of $P_\cF(\check X)$ embeds into
an object of $Inj$. 
 
  We showed in the previous proof that $\nabla^\udot _\alpha\la k\ra =
  j_{\alpha *}\Theta_\alpha[\dim O_\alpha]\la k\ra$ embeds into
  $I^\udot_\alpha\la k\ra$.  
  Let $E_\alpha$ be a $\cF$-local system on $O_\alpha$.  It
  embeds into a finite direct sum of injective $\cF$-local systems
  $\Theta_\alpha\la n\ra$, so $j_{\alpha *}E_\alpha[\dim O_\alpha]$
  embeds into a finite sum of $I^\udot_\alpha\la k\ra$.
  
  We prove that an arbitrary $S^\udot \in P_\cF(\check X)$ embeds into
  an object of $Inj$ by induction on the number of orbits in $\supp
  S^\udot$.  The case when $\supp S^\udot$ is a single orbit 
follows from the previous paragraph.  Otherwise, 
take an open orbit $O_\beta$ in $\supp S^\udot$, and let
  $\phi\colon S^\udot \to \bR j_{\beta *}(S^\udot|_{O_{\beta}})$ be
  the adjunction map.  We have seen that the target of $\phi$ 
embeds into an object $I^\udot_1 \in Inj$, so the image of $\phi$
does as well.  Since $\ker \phi$ has support strictly smaller than
$S^\udot$, by induction it embeds into an object $I^\udot_2\in Inj$.
Since $I^\udot_2$ is injective, this embedding extends to a map
$S^\udot \to I^\udot_2$.  Thus we get an embedding of $S^\udot$
into $I^\udot_1 \oplus I^\udot_2 \in Inj$.
  
The argument for $P_{cf}(\check X)$ is essentially the same.
\end{proof}

Now define
 \[R = \eend(I^\udot)^{opp} \cong
   \oplus_{n\in \Z}\Hom_{\cA\mod}(\cL,\cL\{n\})^{opp}.\]
By Lemma \ref{Noetherian}, $R$ is (left and right) Noetherian.
Then applying Propositions
\ref{mixed categories are modules} and \ref{ungraded modules}, we
conclude that $P_\cF(\check X)$ is equivalent to $R\mcf$,
$P_{cf}(\check X)$ is equivalent to $R\Mcf$, and $F_{cf}$ is the
functor of forgetting the grading.  

\begin{cor} \label{unipotent equivalence}
  There are equivalences of triangulated categories:
  $D^b(LC_\cF(\check X))\simeq D^b(P_\cF(\check X))$ and
  $D^b(LC_{cf}(\check X))\simeq D^b(P_{cf}(\check X))$.
\end{cor} 
The argument for the two equivalences is the same, so we concentrate
on the first one.  Both $D^b(LC_\cF(\check X))$ and $D^b(P_\cF(\check
X))$ are generated by the injectives in $P_\cF(\check X)$; note that
any complex has a bounded injective resolution, since any object of
$D^b(\cA)$ can be represented by a bounded complex of pure sheaves.
We thus need to show that for any injectives $I^\udot_1$, $I^\udot_2
\in P_\cF(\check X)$
and any $d \in \Z$ there is an isomorphism
\[\Ext^d_{P_\cF(\check X)}(I^\udot_1, I^\udot_2) \stackrel{\sim}{\to} 
\Hom_{D^b(LC_\cF(\check X))}(I^\udot_1, I^\udot_2[d]).\] Both sides
are automatically isomorphic for $d \le 0$, while for $d > 0$ the left
side vanishes by the injectivity of $I^\udot_2$.  The vanishing of the
right side is just (*) from the proof of Theorem \ref{K(simple) is
  injective}.

\subsection{Proof of Theorem \ref{K is Koszul}, part II: 
$K^{-1}(\text{simple})$ is projective}
For any $\tau \in [\sig]$, let $\cP^\udot_\tau = K^{-1}(L^\udot_{\tau^\bot})$.  
Analogously to Theorem \ref{K(simple) is injective},
we have
\begin{thm} \label{K(projective) is simple}
$\cP^\udot_\tau$ lies in the core of the perverse $t$-structure on 
$D^b(\cA)$ defined in \S\ref{t-structure on cA-complexes}.  
In the abelian category $P(\cA)$, it is the projective cover
of $\cL^\tau$. 
\end{thm}

\begin{proof}  Let $L^\udot = L^\udot_{\tau^\bot}$. 
Consider an injective resolution of $L^\udot$ (as remarked before, it
can be chosen to be bounded):
\[L^\udot \stackrel\sim\to (J^\udot_0 \to J^\udot_1 \to \dots \to J^\udot_k).\]
Taking $K^{-1}$ gives a complex $\cM_0\to \cM_1 \to\dots\to \cM_k$ of pure 
$\cA$-modules which represents the object $\cP^\udot_\tau$.

This complex will be perverse if $\cM_j[-j]$ is perverse
for $j = 0,\dots,k$, or in other words, if each $J^\udot_l$ is a direct sum 
of objects $I^\udot_{\alpha}\la l\ra$, $\alpha \in \check \Sig$.  
The existence such a resolution follows from Proposition 
\ref{unipotent Koszul} and Corollary \ref{unipotent equivalence}.

Since objects of $P(\cA)$ have finite length, to show that $\cP^\udot_\tau$ 
is the projective cover of $\cL^\tau$ it will be enough to show that for any 
$\rho \in \Sig$, $k,l \in \Z$ we have
$\Hom_{D^b(\cA)}(\cP^\udot_\tau, \cL^\rho[k]\la l\ra)$ is one-dimensional
if $k = l = 0$ and $\rho = \tau$, and vanishes otherwise.  
By applying $K$, this follows from Theorem \ref{K(simple) is injective}.
\end{proof}

Define  $\cP^\udot := \oplus_{\tau\in[\sig]} \cP^\udot_\tau$, and let 
\[\check R = \eend_{P(\cA)}(\cP^\udot)^{opp}.\]

\begin{cor} $\cP^\udot$ is a graded projective
generator of $P(\cA)$; $F_T(\cP^\udot)$ is a projective generator of $P_T(X)$.

There are equivalences of abelian categories $P(\cA) \simeq \check R\mof$ and
$P_T(X) \simeq \check R\Mof$; with respect to these equivalences $F_T$ is
the functor of forgetting the grading. 
\end{cor}

\begin{cor} There are equivalences of triangulated categories: $D^b(\cA)\simeq D^b(P(\cA))$
and $D^b_T(X) \simeq D^b(P_T(X))$.
\end{cor}

The proofs are the same as in the previous section; note that since objects of
$P(\cA)$ have finite length the ring $\check R$ is automatically Noetherian.

\section{Some proofs} \label{appendix}
\subsection{Proof of Theorem \ref{purity of simples}}
The functor $[k]\la -k\ra$ on $D^b(LC_\cF(X))$ preserves the 
property of being pure of weight $0$, so we can instead prove that
\[S^\udot := L^\udot_\alpha[-c(\alpha)]\la c(\alpha)\ra
= j_{\alpha !*}\C_\alpha\] is pure of weight $0$.  The support of
$S^\udot$ is $\ol{O_\alpha}$, which is itself a toric variety (for a
smaller torus).  Thus we can assume that $\alpha = o$ is the zero cone
and the support of $S^\udot$ is all of $X$.

Let $\bS^\udot = F_{cf}S^\udot$; it is an
intersection cohomology sheaf shifted so that the restriction
to the open orbit $O_o$ is a local
system in degree $0$.  The $\cF$-structure on $S^\udot$
defines an isomorphism $\theta\colon \cF^{-1}\bS^\udot
\stackrel{\sim}\to \bS^\udot$.  It induces an action on the stalk of
the cohomology sheaves $H^i(j_\beta^*\bS^\udot)$ and
$H^i(j_\beta^!\bS^\udot)$; we need to show this action is
multiplication by $2^{i/2}$.

Note that if $O_\beta$ has positive dimension, there is an
$\cF$-stable normal slice to $O_\beta$ at a point of $(O_\beta)^\cF$
which is itself an affine toric variety. By restricting to this slice
we can restrict to the case when $O_\beta = \{b\}$ is a single point.

Note that $\cF^{-1}\bS^\udot \cong \bS^\udot$ (see \cite{BM}), and all
automorphisms of $\bS^\udot$ are multiplication by scalars.  Therefore
$\theta$ is uniquely determined by its action on the stalk at an
$\cF$-fixed point of the open orbit $O_o$, where it acts as the
identity.

Let $\pi \colon \wt{X} \to X$ be a toric resolution of singularities,
and let $\wt{\cF}$ be our geometric Frobenius map on $\wt{X}$; we have
$\wt\cF\pi = \pi\cF$.
 
Since $\wt\cF^*\C_{\wt X} \cong \C_{\wt X}$, we can put a
$\wt\cF$-structure on the constant sheaf $\C_{\wt{X}}$ by letting
$\tilde\theta\colon \wt{\cF}^*\C_{\wt{X}} \to \C_{\wt{X}}$ act as the
identity on the stalk at a point of $(O_o)^\cF$.  By adjunction
$\tilde \theta$ induces a map $\C_{\wt{X}} \to
\bR\wt{\cF}_*\C_{\wt{X}}$, and applying $\bR\pi_*$ and adjunction
again gives 
$\theta'\colon \cF^*\bR\pi_*\C_{\wt X}\to \bR\pi_*\C_{\wt X}$.

By the decomposition theorem \cite{BBD}, $\bS^\udot$ is a direct
summand of $\bR\pi_*\C_{\wt X}$ and of $\cF^*\bR\pi_*\C_{\wt X}$.
Composing $\theta'$ with the inclusion and projection gives a map
$S^\udot \to S^\udot$; it is easy to see that it agrees with $\theta$
on the open orbit, so it must equal $\theta$ on all of $X$.

The cohomology groups of $j_\beta^*\bR\pi_*\C_{\wt X}$ and
$j_\beta^!\bR\pi_*\C_{\wt X}$ are $H^\udot(\pi^{-1}(b))$ and
$H^\udot(\wt{X},\wt{X}\setminus \pi^{-1}(b))$, respectively, and the
action of $\theta'$ is the action of the pullback $\wt{\cF}^*$.  Thus
we have reduced the proof of the theorem to showing that this action
is multiplication by $2^{i/2}$ on the cohomology in degree $i$.

Here is one way to see this: $\wt{X}$ has a completion to a smooth
complete toric variety $Y$.  There is a homomorphism $\C^* \to T$ so
that the induced action of $\C^*$ on $X$ is ``attractive'':
$\lim_{t\to 0} t\cdot x = b$ for all $x\in X$, and the induced action
on $Y$ has isolated fixed points.  Then by Bia\l ynicki-Birula
$\pi^{-1}(b)$ has a decomposition into $\bigcup_x {C_x}$ into affine
cells, so $H^\udot(\pi^{-1}(b)) \cong \oplus_x H^\udot_c(C_x)$.  The
cells are $T$- and $\cF$-invariant, and $\wt\cF$ acts on each
$k$-dimensional cell as the map $\C^k \to \C^k$,
$(x_1,\dots,x_k)\mapsto (x_1^2,\dots,x_k^2)$.  The result for
$H^\udot(\pi^{-1}(b))$ follows immediately.

For $H^\udot(\wt{X},\wt{X}\setminus \pi^{-1}(b))$, we use the Bia\l
ynicki-Birula cells for the opposite character $\C^* \to T$.  Then
$\wt{X}$ is an open union of these cells which deformation retracts
onto $\pi^{-1}(b)$ by our action.  Therefore we have
$H^\udot(\wt{X},\wt{X}\setminus \pi^{-1}(b)) \cong H^\udot_c(\wt{X})$, and
can use the argument of the previous paragraph.

\subsection{Proof of Theorem \ref{weight filtration}}
We begin by defining the filtration $W_\udot S^\udot$ 
when $S^\udot \in P_{\cF,c}(X)$, i.e.\
when $S^\udot$ has finite length.  We proceed by induction on the length
of $S^\udot$.  If $S^\udot$ has length $1$, it is simple, say of weight 
$m$, and we can let $W_k S^\udot = 0$ if $k < m$, $W_k S^\udot = S^\udot$
for $k \ge m$. 

Otherwise suppose the filtration has already been defined for objects of 
smaller length.  
Find a simple subobject $L^\udot$ of $S^\udot$,
and suppose it is pure of weight $m$.  Let $\phi \colon 
S^\udot \to C^\udot = S^\udot/L^\udot$ be the corresponding 
quotient map. By induction we can assume 
we have already defined our filtration on $C^\udot$.

For any $k < m$ consider the 
exact sequence
\[0 \to L^\udot \to \phi^{-1} W_k C^\udot \to W_k C^\udot \to 0.\]
Since the simple constituents of $W_k C^\udot$ are all pure of weights
$< m$, Proposition \ref{unipotent Koszul} implies
\[\Hom(L^\udot,  W_k C^\udot) = \Ext^1(L^\udot,  W_k C^\udot) = 0,\]
and so the exact sequence splits canonically.  We then define
$W_kS^\udot$ to be the image of $W_kC^\udot \to  \phi^{-1} W_k C^\udot
\to S^\udot$, where the first map is the splitting map.
For $k \ge m$ we let $W_k S^\udot = \phi^{-1}(W_k C^\udot)$.
Then $\Gr^W_m S^\udot \cong \Gr^W_m C^\udot \oplus L^\udot$, while
$\Gr^W_k S^\udot \cong \Gr^W_k C^\udot$ if $k \ne m$, so 
 $\Gr^W_k S^\udot$ is pure of weight $k$ for all $k\in\Z$, since
the same was true for $C^\udot$ by induction.  

Next we extend this filtration to arbitrary objects of $P_\cF(X)$,
which may not have finite length.  In order to do this, we need to
show that for any object $S^\udot\in P_\cF(X)$ there is a lower bound
on the weights of the simple constituents of $S^\udot$.  To see this,
use induction on the number of orbits in the support of $S^\udot$.  If
the support is a single orbit $O_\alpha$, the result follows from the
equivalence $P_\cF(O_\alpha) \simeq LC_\cF(O_\alpha) \simeq
\co\cT_\alpha\mcf$.

Otherwise, let $O_\alpha$ be an open orbit in the support of
$S^\udot$.  Let $\phi$ denote the natural adjunction morphism $S^\udot
\to \bR j_{\alpha*}(S^\udot|_{O_\alpha})$, and consider the short
exact sequence
\[0 \to \ker \phi \to S^\udot \to \im \phi \to 0.\]
The support of $\ker \phi$ is strictly smaller, so its weights are
bounded below by the inductive hypothesis.  Thus it will suffice to
show the weights of $\im \phi$ are bounded below as well.  But by the
preceding paragraph the weights in $S^\udot|_{O_\alpha}$ are bounded
below, say by $w$.  Since $\bR j_{\alpha *}$ is a $t$-exact functor,
this implies that the weights of $\bR
j_{\alpha*}(S^\udot|_{O_\alpha})$, and hence of $\im \phi$, are
bounded below by $w+w'$, where $w'$ is a lower bound for the weights
of $\bR j_{\alpha*}\C_\alpha[\dim O_\alpha]$.  This lower bound exists
because $\bR j_{\alpha*}\C_\alpha[\dim O_\alpha]$ has finite length by
Proposition \ref{constructible = finite length} (in fact, $w' = 0$
works).

The existence of the filtration $W_\udot$ follows immediately: if
$S^\udot \in P_\cF(X)$ and $k\in\Z$, the collection $W_k\hat
S^\udot$ forms a directed system over all finite-length subobjects
$\hat S^\udot$ contained in $S^\udot$.  It vanishes identically for 
$k \ll 0$, so we can proceed by induction on $k$: assume that
$W_{k-1}S^\udot$ has been defined in such a way that $S_+^\udot =
S^\udot/W_{k-1}S^\udot$ has only simple constituents of weights $\ge
k$.  The family $\{W_kS^\udot_+\}$ must be eventually constant, since
$\Hom(L^\udot_\alpha\la k\ra, S^\udot_+)$ is finite-dimensional for
all $\alpha$.  Thus $\{W_kS^\udot\}$ also stabilizes, so the limit is
a finite-length subobject.  

The properties of a mixed category are easy to verify.

\subsection{Proof of Theorem \ref{t-exact}} 
  Let us briefly recall how the functor $\epsilon\colon
  D^b_T(X) \to D^b(\cA_X\mof)$ of \cite{L} is defined.  Since $X$
  is affine, we can choose a $T$-equivariant embedding $X \hookrightarrow
  \PP^n$, where the action of $T$ on $\PP^n$ is linear.  We can choose
  a representative for the classifying space $BT$ so that $\PP^n_T$ is
  an infinite dimensional manifold in the sense of \cite{BL} --
  essentially this means it is a limit of finite-dimensional manifolds
  by closed embeddings.  $\PP^n_T$ has a ``de Rham complex''
  $\Omega^\udot_{\PP^n_T}$ which is a resolution of the constant sheaf
  $\R_{\PP^n_T}$ by soft sheaves.  It is also naturally a
  supercommutative sheaf of DG-algebras.  We then let
  $\Omega^\udot_{X_T} = \Omega^\udot_{\PP^n_T}|_{X_T}.$
  
  Let $\pi\colon X_T\to X/T$ be the map sending $O_T$ to $O/T$ for any
  $T$-orbit $O$.  Given $S^\udot \in D^b_T(X)$, %, which we think of as a
  %complex of injective sheaves on $X_T$, 
  the complex $M^\udot = \pi_*(\Omega^\udot_{X_T} \otimes S^\udot)$ is
  naturally a DG-module over the DG-sheaf $\wt\cA:=
  \pi_*(\Omega^\udot_{X_T})$.  Here $\otimes$ is tensoring over $\R$.
  Since all $\R$ sheaves are flat, $\Omega^\udot_{X_T} \otimes
  S^\udot$ is quasi-isomorphic to $S^\udot$.
  
  The DG-sheaf $\wt\cA$ is formal, i.e.\ it is quasi-isomorphic to its
  cohomology $H(\wt\cA)$.  Under the natural identification of $X/T$
  with the fan $\Sig$ defining $X$, $H(\wt\cA)$ is canonically
  isomorphic to our sheaf of rings $\cA_\Sig$.  This gives an
  equivalence of categories
\begin{equation}\label{xxx}
D(\text{DG-$\wt\cA$})\simeq D(\text{DG-$\cA$})
\end{equation}
which commutes with restriction and corestriction.  The functor
$\epsilon$ is the composition of this equivalence with
$\pi_*(\Omega^\udot_{X_T} \otimes \udot)$.

\newcommand{\hatj}{\hat\jmath}

Let us prove (a).  Since the functor $F_{T,\Sig}$ is defined locally,
we can assume that $\Sig = [\sig]$, so $O_\sig$ is the unique closed
orbit in $X$.  Let $S^\udot = F_{T,\Sig}\cM^\udot$.  Let $j\colon
O_{\sig,T} \to X_T$, $\hatj\colon \{\sig\}\to \Sig$ denote the
inclusions, and let $\pi_\sig$ be the constant map $O_\sig \to
\{\sig\}$, so $\pi\circ j = \hatj\circ \pi_\sig$. We will show that
there are quasi-isomorphisms
\[\hatj^*\pi_*(\Omega^\udot_{X_T} \otimes S^\udot) \simeq 
\pi_{\sig *}j^*(\Omega^\udot_{X_T} \otimes S^\udot) \simeq \pi_{\sig
  *}(\Omega^\udot_{O_{\sig,T}}\otimes j^*S^\udot).\] This will imply
our result -- the equivalence \eqref{xxx} commutes with taking stalks,
so the left hand side is $\hatj^*\nu\cM^\udot = \nu(\cM(\sig))$, while
the right hand side is $\epsilon(j^*F_{T,\Sig}\cM^\udot)$.  Applying
$\epsilon^{-1}$ gives (a).

The second isomorphism is standard; see \cite[Proposition 2.3.5]{KS},
for instance.  For the first isomorphism, note that since the smallest
open subset of the fan $[\sig]$ containing $\sig$ is $[\sig]$ itself,
the functor $\hatj^*$ is naturally isomorphic to $\hat p_*$, where
$\hat p\colon [\sig]\to \{\sig\}$ is the constant map.  Therefore it
will be enough to construct a quasi-isomorphism
\[\hat p_*\pi_*(\Omega^\udot_{X_T} \otimes S^\udot) = 
\pi_{\sig *}p_*(\Omega^\udot_{X_T} \otimes S^\udot)
\stackrel{\sim}{\to} \pi_{\sig *}j^*(\Omega^\udot_{X_T} \otimes
S^\udot),\] where
\[p = p_{\sig,T} \colon X_T \to O_{\sig,T}\] is the map induced by
the projection map $p_\sig$ defined in \S\ref{toric projections}.

Since $p\circ j$ is the identity on $O_{\sig,T}$, adjunction gives a 
natural transformation $p_* \to p_*j_*j^* = j^*$.  We will show that applying
it to $\wt S^\udot = \Omega^\udot_{X_T} \otimes S^\udot$ gives a
quasi-isomorphism.  Without loss of generality we can assume that
$S^\udot = \bR i_{\tau *}\R_{O_{\tau,T}}$, where $\tau \in [\sig]$ and
$i_\tau\colon O_{\tau,T} \to X_T$ is the inclusion.  Note that $\wt
S^\udot$ is a complex of soft sheaves, so applying $p_*$ to it is the
same as applying $\bR p_*$.  The stalk cohomology of $\bR p_*\wt S^\udot$
and $j^*\wt S^\udot$ at a point $x\in O_{\sig,T}$ are both isomorphic
to the cohomology of the torus $p^{-1}(x) \cap O_{\tau,T}$, which
implies the claim.

For (b), consider the chain of maps
\[\bR \pi_{\sig *}(j^*\Omega^\udot_{X_T} \otimes j^!S^\udot) \to
\bR \pi_{\sig *}j^!(\Omega^\udot_{X_T}\otimes S^\udot)
\stackrel{\sim}{\to} \hatj^! \bR\pi_*(\Omega^\udot_{X_T}\otimes
S^\udot).\] For the first map see \cite[Proposition 3.1.11]{KS}; the
second map is the usual base change.  To check this is an isomorphism
it is enough to consider the case $S^\udot = \bR i_{\tau *}
\R_{O_{\tau,T}}$ as before.  If $\tau = \sig$ the isomorphism is
clear, while if $\tau \ne \sig$ both sides vanish.  (b) then follows.

The isomorphisms (c) and (d) follow from these statements using
results of \cite{BL}.  In the case of a single orbit $O_\sig$,
the equivalence 
$\epsilon\colon D^b_T(O_\sig) \to D_f(\text{DG-}\cA_\sig)$
can be factored as an equivalence 
\begin{equation} \label{yyy}
D^b_T(O_\sig) 
\stackrel\sim\to D^b_{T/T_\sig}(pt)
\end{equation}
\cite[Theorem 2.6.2]{BL} (here $T_\sig$ is the stabilizer of any
point of $O_\sig$) followed by an equivalence 
$D^b_{T/T_\sig}(pt) \stackrel\sim\to D_f(\text{DG-}\cA_\sig)$
\cite[Theorem 12.7.2(ii)]{BL}.  The pullback functor $i^*_x$
is $Q^*_f$, where $f\colon \{y\} \to O_\sig$ 
is the inclusion of a point $\{y\}$, which is a $\phi$-map for
the homomorphism $\phi\colon \{1\} \to T$ (for the definition and
properties of $Q^*_f$, see \cite[\S3.6]{BL}).  The equivalence
\eqref{yyy} is $Q^*_g$, where $g\colon O_\sig\to pt$ is the quotient
map and $pt$ carries a $T/T_\sig$-action.  Thus $i^*_x = Q^*_fQ^*_g = 
Q^*_{gf}$.  Applying \cite[Theorem 12.7.2(iii)]{BL} completes the proof.
\bibliography{eqkd}
\end{document}